\documentclass[11pt,reqno]{amsart}
\usepackage{amssymb,mathrsfs,graphicx,subfigure, enumerate}
\usepackage{amsmath,amsfonts,amssymb,amscd,amsthm,bbm}
\usepackage{graphicx,colortbl,here}
\usepackage{extpfeil}

\def\cD {\mathcal{D}}
\def\cE {\mathcal{E}}

\def\cI {\mathcal{I}}

\def\cK {\mathcal{K}}
\def\cL {\mathcal{L}}
\def\cM {\mathcal{M}}

\def\cS {\mathcal{S}}

\def\ka {{\kappa}}

\def\d {{\partial}}

\def \rpsi_i {|\psi_i \rangle}
\def \lpsi_i {\langle \psi_i|}
\def \lrpsi_i{\langle \psi_i | \psi_i \rangle}
\def \rpsi_k {|\psi_k \rangle}
\def \lpsi_k {\langle \psi_k|}
\def \lrpsi_k{\langle \psi_k | \psi_k \rangle}

\newcommand{\bbr}{\mathbb R}

\newcommand{\bbs}{\mathbb S}
\newcommand{\bbz}{\mathbb Z}
\newcommand{\bbh}{\mathbb H}

\newcommand{\kp}{\kappa}

\newcommand{\St}{\textup{St}(p,n)}
\newcommand{\tF}{\textup{F}}

\newcommand{\veps}{\varepsilon}

\newcommand{\ba}{\begin{aligned}}
\newcommand{\ea}{\end{aligned}}

\newcommand{\be}{\begin{equation}}
\newcommand{\ee}{\end{equation}}

\newtheorem{theorem}{Theorem}[section]
\newtheorem{lemma}{Lemma}[section]
\newtheorem{corollary}{Corollary}[section]
\newtheorem{proposition}{Proposition}[section]
\newtheorem{remark}{Remark}[section]
\newtheorem{definition}{Definition}[section]

\renewcommand{\d}{\textup{d}}

\begin{document}

\title[Emergent dynamics on Stiefel manifolds]{Emergent behaviors of high-dimensional Kuramoto models on Stiefel manifolds}

\author[S-Y. Ha]{Seung-Yeal Ha}
\address[S.-Y. Ha]{\newline Department of Mathematical Sciences and Research Institute of Mathematics \newline Seoul National University, Seoul 08826  and \newline
Korea Institute for Advanced Study, Hoegiro 85, Seoul, 02455, Republic of 
Korea} \email{syha@snu.ac.kr}

\author[M. Kang]{Myeongju Kang}
\address[M. Kang]{\newline Department of Mathematical Sciences and Research Institute of Mathematics \newline Seoul National University, Seoul 08826, Republic of Korea}
\email{bear0117@snu.ac.kr}

\author[D. Kim]{Dohyun Kim}
\address[D. Kim]{\newline Department of Mathematics, \newline Sungshin Women's University, Seoul 02844, Republic of Korea}
\email{dohyunkim@sungshin.ac.kr}

\thanks{\textbf{Acknowledgment.}
The work of S.-Y. Ha was supported by NRF grant (2017R1A2B2001864).}

\begin{abstract}
We study emergent asymptotic dynamics for the first and second-order high-dimensional Kuramoto models on Stiefel manifolds which extend the previous consensus models on Riemannian manifolds including several matrix Lie groups. For the first-order consensus model on the Stiefel manifold proposed in [Markdahl et al, 2018], we show that the homogeneous ensemble relaxes the complete consensus state exponentially fast. On the other hand for a heterogeneous ensemble, we provide a sufficient condition leading to the phase-locked state in which relative distances between two states converge to definite values in a large coupling strength regime. We also propose a second-order extension of the first-order one by adding an inertial effect, and study emergent behaviors using Lyapunov functionals such as an energy functional and an averaged distance functional.

%
%
%

\end{abstract}



\date{\today}

\maketitle



\section{Introduction} \label{sec:1}
\setcounter{equation}{0}
Collective behaviors of biological and chemical oscillators have been widely studied not only in applied mathematics but also in other scientific disciplines, for instance, flocking of drones  \cite{H-G, P-E-G} and passivity-based distributed optimization  \cite{Chopra-S, H-C-I-L,  S-S-C-G} in control theory, social dynamics \cite{B-D, B-G, K-W, N} and swarming behavior in quantitative biology \cite{B-B, C-U, D-M-T-H, G-F-M-S, R-C-N}. Despite of its crucial roles in biological processes, the mathematical study of such collective dynamics has been started only after the seminal work of Winfree \cite{Wi2} and Kuramoto \cite{Ku1,Ku2}  a half century ago. Among other models describing collective oscillatory behavior, to name a few, the Cucker-Smale model \cite{C-S}, the Kuramoto model \cite{Chopra-S1, Ku1, Ku2} and the Vicsek model \cite{Vi}, our main interest lies in consensus models on Riemannian manifolds. So far, there has been much available literature dealing  with emergent dynamics on Riemannian manifolds, for instance,  on the unit sphere $\bbs^{d-1}$ in \cite{O1, Z-Z}, on the hyperboloid $\bbh^{d-1}$ in  \cite{R-L-W} and on the matrix Lie groups including special orthogonal group $\textup{SO}(d)$ in  \cite{D-F-M1, S-S-L, T-A-V}, the unitary group $\textup{U}(d)$ in \cite{H-K-R3, H-R, Lo}, the \textit{Lohe group} in \cite{H-K-R2} and  even for the quaternions $\bbh_1$ in \cite{D-F-M-T}. In particular, we are concerned with the Stiefel manifold \cite{Stiefel} in $\bbr^{n \times p}$ for $p\leq n$ consisting of all rectangular matrices satisfying the relations:
\begin{equation*}
\St:= \{ S\in M_{n,p}(\bbr): S^\top  S = I_p\},\quad \|S\|_\text{F}^2 := \textup{tr}(S^\top S) = p,
\end{equation*}
where $S^\top $ is defined as the transpose of a  matrix  $S$ and $M_{n,p}(\bbr)$ denotes the set of all $n\times p$ matrices with real entries. 

Note that the unit sphere and (special) orthogonal group can be recovered from the specific choices of $(p,n)$. Thus, the Stiefel manifold is a  manifold including $\bbs^d$, $\textup{SO}(d)$ and  $\textup{O}(d)$ as special cases. We will briefly review basic properties of the Stiefel manifold in Section \ref{sec:2.1}.  It is worthwhile to mention that an optimization problem $\displaystyle \min_{S \in \St} f(S)$ for an  objective function $f:\St \to \bbr$ has been extensively studied due to its computational difficulty and broad applications, e.g., the linear eigenvalue problem \cite{G-V,W-Y-L-Z}, the nearest low-rank correlation matrix problem \cite{L-Q}, singular value decomposition \cite{L-W-Z,S-I}, Riemmanian optimization \cite{Zhu} and applications to computer vision \cite{Lui,T-V-S-C}. We refer the reader to \cite{E-A-S,S-S,W-Y} and reference therein for introductions to optimization problem on the Stiefel manifold and applications. 

In this paper, we mainly consider the consensus model in \cite{Ma1, Ma2}  in which particles interact with neighboring ones on the Stiefel manifold and for the consistency of the paper, we briefly introduce the model below. Let $(\St, \|\cdot\|)$ be the Stiefel manifold canonically embedded into the Euclidean space $\bbr^{n\times p}$ with its Frobenius (or Euclidean) norm and $\mathcal S:= (S_1,\cdots,S_N)$ denote state ensemble of particles. Next, we define the potential function $\mathcal V$ as the total misfit functional with a symmetric connectivity matrix $\mathcal A := (a_{ik})$: 
\begin{equation*}
\mathcal V(\mathcal S):= \frac{\kp}{N^2} \sum_{i,k=1}^N a_{ik} \|S_i - S_k\|_\tF^2,\quad a_{ik} = a_{ki}>0.
\end{equation*}
Then, the gradient flow with the potential function $\mathcal V$ reads as
\begin{equation} \label{A-0}
\dot S_i  = -\nabla_{S_i} \mathcal V =\frac\kp N  \sum_{k=1}^N a_{ik} \left[ S_k - \frac12 (S_i S_i^\top  S_k + S_i S_k^\top  S_i) \right],\quad i=1,\cdots,N.
\end{equation}
See Section \ref{sec:2.2} for detailed description. We finally add generalized natural frequencies $\Xi_i \in \mathfrak{so}(p)$ into \eqref{A-0} to find the desired model so that the Stiefel manifold is still positively invariant: 
\begin{equation} \label{A-0-0}
\dot S_i  =  S_i \Xi_i + \frac\kp N  \sum_{k=1}^N a_{ik} \left[ S_k - \frac12 (S_i S_i^\top  S_k + S_i S_k^\top  S_i) \right],\quad i=1,\cdots,N.
\end{equation} 
Since the Stiefel manifold is compact, the global existence of a unique solution to \eqref{A-0-0} directly follows from the standard Cauchy-Lipschitz theory. Detailed description will be presented in Section \ref{sec:2.1}.  

For the perspective of optimization, a gradient flow \eqref{A-0} with a total distance as a potential can be regarded an optimization problem of a given target (or objective) function defined on the Stiefel manifold. Thus, our model with a consensus estimate provides a method which tackles an optimization problem on the Stiefel manifold by means of dynamical systems approach.  It should be mentioned that such optimization problems  on manifolds have   extensively studied in literature \cite{B-C-S,  C-M-S-Z, L-B,P-D-M,  S-S}. In particular, consensus-based optimization (CBO for brevity) algorithm toward a global optimization \cite{C-C-T-T, F-H-P-S, F-H-P-S2, G-P, P-T-T-M} has been recently proposed (see \cite{K-K-K-H-Y} for CBO on the Stiefel manifold). On the other hand for a control (or consensus protocol) perspective, it is mentioned in \cite{Ma20a} that \textit{ad-hoc} control algorithm would be employed in a specific situation. In this manner, one can consider \eqref{A-0} as a control problem by regarding the network structure $a_{ik}$  as an external control (or parameter) to obtain a desired pattern formation (see \cite{L-S} on the unit sphere and  \textit{an Olympic ring} in \cite{C-K-P-P} for a distributed control approach for the Cucker-Smale flocking model \cite{C-S}). For instance, if we set $a_{ik}$ to be negative (or competitive), one can expect deployment where all agents tend to splay state, whereas complete consensus would be achieved for a positive (or cooperative) $a_{ik}$. However in this paper, we restrict ourselves  to a cooperative network so that emergence of complete consensus is our primary concern.

We analyze system  \eqref{A-0-0} with all-to-all network topology and left-translation invariance.  More precisely, we choose  all-to-all network $a_{ij} = \kp/N$ where $\kp$ denotes the coupling strength so that all particles communicate with neighbors through the same weight. In addition, we do not consider the effect of  $\Omega_i$ by setting $\Omega_i \equiv O$ due to the left-translation property (see Lemma \ref{leftinv} for details). Then, our first-order consensus model is governed by the following Cauchy problem:
\begin{equation} \label{A-1}
\begin{cases}
\displaystyle \dot S_i =  S_i \Xi_i + \frac{\kp}{N}\sum_{k=1}^N \left[  S_k - \frac12 ( S_iS_i^\top S_k + S_i S_k^\top S_i)\right], \quad t>0,\\
\displaystyle S_i(0) = S_i^\textup{in} \in \text{St}(p,n),\quad i=1,\cdots,N.
\end{cases}
\end{equation}
Although the gradient flow structure guarantees that the Stiefel manifold is positively invariant along the flow \eqref{A-0-0}, we provide its alternative and direct proof in Lemma \ref{L2.1}. 

Next, we turn to the second-order extension of the first-order model \eqref{A-1}. In \cite{H-K2}, the (first-order) Lohe matrix model on the unitary group has been extended to a second-order one by incorporating inertial effect. Similarly, we can also extend  the first-order model \eqref{A-1} into a second-order one whose dynamics is governed by the following Cauchy problem:
\begin{equation} \label{A-2}
\begin{cases}
\displaystyle m  \ddot S_i = - mS_i  \dot S_i^\top  \dot S_i - \gamma  \dot S_i + S_i \Xi_i + \frac{m}{\gamma} ( 2\dot S_i \Xi_i - S_i \Xi_i S_i^\top  \dot S_i + S_i \dot S_i^\top  S_i \Xi_i) \\
\displaystyle \hspace{1.2cm}+ \frac{\kp}{N}\sum_{k=1}^N \left(  S_k - \frac12 ( S_iS_i^\top S_k + S_i S_k^\top S_i)\right),\quad t>0, \\
\displaystyle S_i(0) = S_i^\textup{in} \in \St,\quad \dot S_i^{t,\textup{in}} S_i^\textup{in} + S_i^{t,\textup{in}} \dot S_i^\textup{in} =  O,\quad i=1,\cdots,N,
\end{cases}
\end{equation}
where $m$ and $\gamma$ represent mass and friction, respectively. Although it seems that the model looks quite complicated, we show the the Stiefel manifold is still positively invariant along   system \eqref{A-2}. Of course, if we turn off the inertial effect, that is, $m=0$, then the first-order model \eqref{A-1} can be recovered from \eqref{A-2} straightforwardly. We also note that the Cauchy problems \eqref{A-1} and \eqref{A-2} have a unique global solution due to the compactness of  Stiefel manifold and standard Cauchy-Lipschitz theory. Here, a global solution is referred as a solution which exists for all time (or globally). In other words, a solution does not blow up in finite time. For more detailed description of \eqref{A-2}, we  refer the reader to Section \ref{sec:2.2}.  Next, we recall several concepts for consensus as follows: 
\begin{definition} \label{D1.1}
Let ${\mathcal S} = (S_1, \cdots, S_N)$ be a global solution to \eqref{A-1} or \eqref{A-2}. 
\begin{enumerate}
\item
System \eqref{A-1} or \eqref{A-2} exhibits complete consensus  if the following convergence holds:
\begin{equation*}
\lim_{t\to\infty} \|S_i(t) - S_j(t)\|_\tF = 0,\quad \textup{for all $i,j=1,\cdots,N$}. 
\end{equation*}
\item
System \eqref{A-1} or \eqref{A-2}  exhibits practical consensus  if the following convergence holds:
\begin{equation*}
\lim_{\kp\to\infty} \limsup_{t\to\infty} \|S_i(t) - S_j(t)\|_\tF = 0,\quad \textup{for all $i,j=1,\cdots,N$}. 
\end{equation*}
\item
A global solution to system \eqref{A-1}  or \eqref{A-2} tends to a phase-locked state if the following relation holds:
\begin{equation*}
\lim_{t\to\infty} S_i^\top  S_j \quad \textup{exists for all $i,j=1,\cdots,N$.}
\end{equation*}
\end{enumerate}
\end{definition}
The main results of this paper deal with the emergent collective behaviors for the first-order model \eqref{A-1} and the second-order model \eqref{A-2}. First, we consider the first-order model \eqref{A-1} with a homogeneous ensemble. The corresponding proof can be found in Section \ref{sec:3.1} 
\begin{theorem} \label{T1.1}
Suppose that initial data and system parameters satisfy 
\begin{equation*}
\mathcal D(\mathcal S^\textup{in})<\sqrt2,\quad \Xi \equiv O,\quad a_{ik}:\textup{undirected and connected graph},
\end{equation*}
and let $\mathcal S$ be a global solution to \eqref{A-0-0}. Then, system \eqref{A-0-0} exhibits complete consensus exponentially fast.
  \end{theorem}

%
The results in \cite{Ma1, Ma2} deal with almost global consensus without any explicit decay estimate. On the other hand, for a heterogeneous ensemble, we establish the emergence of the locked state exponentially fast in a large coupling regime. 
\begin{theorem} \label{T1.2}
Suppose that the coupling strength and initial data satisfy \eqref{C-16} and let $\mathcal S$  be a  global solution to \eqref{A-0-0}. Then, $\mathcal S$ tends to a locked state.
\end{theorem}
For a detailed initial framework and proof, we refer to Section \ref{sec:3.2}.\newline

%

Next, we turn to the second-order model \eqref{A-2}. As in the first-order one, we consider both homogeneous and heterogeneous ensembles. For the desired results, we first derive an energy estimate (Proposition \ref{P4.1}):
\begin{align} \label{A-3}
\begin{aligned}
& \frac{\d}{\d t}\left( \frac{m}{N} \sum_{i=1}^N \|\dot{S}_i \|^2_\tF +\frac{\kappa}{2N^2} \sum_{i,j=1}^N \left\|S_i-S_j\right\|_\tF^2\right) \\
& \hspace{1cm} = -\frac{2\gamma}{N}\sum_{i=1}^N\| \dot{S}_i \|_\tF^2 +\frac{1}{N}\sum_{i=1}^N\textup{tr} ( \dot{S}_i^\top S_i\Xi_i -\Xi_i S_i^\top \dot{S}_i),\quad t>0. 
\end{aligned}
\end{align}
In what follows, we assume that the network topology satisfies
\begin{align*}
&a_m := \min a_{ik},\quad a_M := \max a_{ik},\quad d(\mathcal A) := \max |a_{ik} - a_{jk}|, \\
&\frac1N \sum_{k=1}^N a_{ik} \equiv \xi,\quad i=1,\cdots,N,\quad  0<  \Lambda :=   a_m - \frac{N-1}{N}(a_M + d(\mathcal A)) < 8pa_M^2.
\end{align*}
Furthermore, for the averaged relative distances $\displaystyle \mathcal G := \frac{1}{N^2} \sum_{i,j=1}^N \|S_i - S_j\|_\tF^2$, we derive a second-order differential inequality (Lemma \ref{L4.2}) for $\mathcal G$:
\begin{equation}\label{A-4}
m \ddot {\mathcal G} + \gamma \dot{\mathcal G} + 2\kp \xi \mathcal G  \leq 16mD(\dot{\cS})^2 +8\|\Xi\|_\infty    +\frac{16m\sqrt{p}\|\Xi\|_\infty}{\gamma}D(\dot{\cS}),
\end{equation}
where $\mathcal D(\dot S)$ and $\|\Xi\|_\infty$ are defined as follows:
\begin{equation*}
\mathcal  D(\dot S):=\max_{1\leq i \leq N } \|\dot S_i\|_\tF,\quad \|\Xi\|_\infty: = \max_{1\leq i \leq N } \| \Xi_i\|_\tF. 
\end{equation*}
Based on two estimates \eqref{A-3} and \eqref{A-4}, for a homogeneous ensemble, we show that the complete consensus occurs for some admissible initial data.
\begin{theorem}\label{T1.3}
Suppose that system parameters and initial data satisfy
\begin{align*}
&m >0,\quad \gamma > 0,\quad \kappa > 0,\quad \Xi_i \equiv O \quad \textup{for}\ i=1,\cdots,N,\\
&\mathcal E(0) = \frac mN\sum_{i=1}^N \|\dot S_i^0\|_\tF^2 + \frac{\kp}{2N^2} \sum_{i,j=1}^Na_{ij}  \|S_i^0 - S_j^0\|_\tF^2 <\infty,
\end{align*}
and let $\cS$ be a global solution to $\eqref{A-2}$. Then, we have
\begin{align*}
\lim_{t\to\infty}\| \dot{S}_i(t) \|_\tF = 0\quad \textup{for}\ i = 1,\cdots,N.
\end{align*}
Moreover, system \eqref{A-2} exhibits the complete consensus:
\begin{align*}
\lim_{t\to\infty} \mathcal G(t)  = 0.
\end{align*}
\end{theorem}
We refer the reader to Section \ref{sec:4.1} for the proof. 
In contrast, for a heterogeneous ensemble, we assume that the inertia and the coupling satisfy the following relation:
\begin{equation*}
m\kp^{1+\eta} = \mathcal O(1) \quad \textup{for some $\eta>0$}. 
\end{equation*} 
Then, under this setting, we arrive at the following result.
\begin{theorem} \label{T1.4}
Suppose that  system parameters and initial data satisfy 
\begin{equation} \label{Y-1}
D(\dot S^\textup{in}) < \frac1\gamma ( \|\Xi\|_\infty + \kp a_M \sqrt p ), \quad m = \frac{m_0}{\kp^{1+\eta}},
\end{equation}
and let $\mathcal S$ be a global solution to \eqref{A-2}. Then, system \eqref{A-2} exhibits practical consensus:
\begin{equation*}
\lim_{\kp\to\infty} \limsup_{t\to\infty} \mathcal G(t) = 0.
\end{equation*}
\end{theorem}
The proof can be found in Section \ref{sec:4.2}. 
 
 \vspace{0.5cm}

 The rest of the paper is organized as follows. In Section \ref{sec:2}, we briefly discuss properties of the Stiefel manifold to be used later and provide descriptions of the first-order and second-order models. In Section \ref{sec:3}, we  present proofs of Theorem \ref{T1.1} and Theorem \ref{T1.2} which deal with the first-order model. In Section \ref{sec:4}, rigorous justification of Theorem \ref{T1.3} and Theorem \ref{T1.4} for the second-order model is provided. Finally, Section \ref{sec:5} is devoted to a brief summary of our main results and discussion for a future work. 

\vspace{0.5cm}

\noindent{\bf Notation:} We denote by $M_{n,p}(\bbr)$ as the set of all $n\times p$ real matrices and for notational simplicity, we set $M_n(\bbr) := M_{n,n}(\bbr)$. In addition, $O$ is the  zero matrix regardless of its size. 


\section{Preliminaries} \label{sec:2}
\setcounter{equation}{0}
In this section, we briefly discuss the Stiefel manifold and present detailed description and properties  of the first-order and second-order consensus models on the Stiefel manifold to be used later in later sections. 
\subsection{The Stiefel manifold} \label{sec:2.1}
We define Stiefel manifold and  Frobenius norm:
\[
\text{St}(p,n) := \{ S\in M_{n,p}(\bbr): S^\top  S = I_p\}, \qquad \|S\|_\text{F}^2 := \text{tr}(S^\top S) = p\quad \textup{for $S \in \St$}. 
\]
Alternatively, it can be defined as the set of all $p$-tuples $(x_1,\cdots,x_p)$ of orthonormal vectors in $\bbr^n$  or it is isomorphic to a homogeneous space:  
\begin{equation*}
\St \simeq \textup{O}(n) / \textup{O}(n-p) .
\end{equation*}
In addition, if $p$ is strictly less than $n$, then one also finds
\begin{equation*}
\St \simeq \textup{SO}(n) / \textup{SO}(n-p).
\end{equation*}
Thus, the Stiefel manifold $\St$ is a compact set whose dimension is $pn - p(p+1)/2$. Furthermore, it is well known that $\St$ reduces to several well-known manifolds, for instance, 
\begin{equation} \label{B-4}
\text{St}(1,n) = \bbs^{n-1} \subseteq \bbr^{n},\quad \text{St}(n-1,n) = \text{SO}(n),\quad \text{St}(n,n) = \text{O}(n). 
\end{equation}
We set $\mathfrak{so}(n)$ to be the special orthogonal Lie algebra associated with $\text{SO}(n)$. Then, we define  two maps $\text{skew}:M_n(\bbr) \to \mathfrak{so}(n)$ and $\text{sym}:M_n(\bbr) \to \mathfrak{so}(n)^\perp$ as
\begin{equation*}
\text{skew}(X) := \frac12 (X- X^\top ),\quad \text{sym}(X) := \frac12 ( X+ X^\top ).
\end{equation*}
Then, it is easy to see that the tangent space and the normal space of $\text{St}(p,n)$ at a point $S$ are defined by 
\begin{align*}
&\textup{T}_S\St := \{ A \in M_{n,p}(\bbr): \text{sym} (S^\top  A ) = O \} = \{ A \in M_{n,p}(\bbr): S^\top   A + A^\top  S = O \}, \\
&\textup{N}_S\St : = \{ SV: \textup{$V$ is a $p\times p$ symmetric matrix}\},
\end{align*}
and the projection of $X$ onto $\textup{N}_S\St$  is given by  $S\textup{sym}(S^\top X)$. Thus, the projection by $\Pi:M_{n,p}(\bbr)\times \St \to \textup{T}_S\St$ is written as
\begin{equation*}
\Pi(X,S) =  X- S \text{sym}(S^\top X) = S \text{skew}(S^\top X) + (I_n - SS^\top )X. 
\end{equation*}
For further details for the Stiefel manifold, we refer to \cite{E-A-S}. 

\subsection{A first-order consensus model on $\text{St}(p,n)$}  \label{sec:2.2} 
In this subsection, we review the first-order model proposed in \cite{Ma1,Ma2} and study its basic property. First, we state the positive invariance of $\St$ for \eqref{A-0-0} which can be guaranteed from the gradient flow structure. 

 \begin{lemma}[Positive invariance of the Stiefel manifold] \label{L2.1}
Let $\mathcal S$ be a global solution to \eqref{A-0-0} with the initial data $\mathcal S^\textup{in} := (S_1^\textup{in}, \cdots, S_N^\textup{in})$. Then, we have 
\begin{equation*}
S_i^\textup{in} \in \St \quad \Longrightarrow  \quad S_i(t) \in \St,\quad t>0.
\end{equation*}
\end{lemma}

Next, we consider the left-translation invariance property whose proof directly follows from straightforward calculations.

\begin{lemma}[Left-translation invariance] \label{leftinv}
For all $L \in \textup{O}(n)$, system \eqref{A-0-0} is invariant under left-translation by an $n \times n$ orthogonal matrix in the sense that a transformed variable $V_i:= LS_i$ satisfies
\begin{equation*}
\dot V_i =   \Omega_i V_i + V_i \Xi_i  + \sum_{k=1}^N a_{ik} \left( V_k - \frac12(V_iV_i^\top V_k + V_iV_k^\top V_i)\right).
\end{equation*}
\end{lemma}

The model \eqref{A-0-0} on the Stiefel manifold in fact includes several first-order models  on the Riemannian manifolds such as  $\bbs^{n-1}$ in \cite{O1}, $\textup{SO}(n)$ in \cite{S-S-L} and $\bbs^1$ in \cite{Ku1}: 
\begin{align} \label{B-7-1}
\begin{aligned}
&\text{(i)}~~ \dot R_i = \Omega_i R_i + \sum_{j=1}^N a_{ij} R_i \text{skew}(R_i^\top R_j),\quad R_i \in \text{SO}(n). \\
&\text{(ii)}~~ \dot x_i = \Omega_i x_i + (I_{n} - x_i \otimes x_i) \sum_{j=1}^N a_{ij} x_j,\quad x_i \in \bbs^{n-1}. \\
&\text{(iii)}~~ \dot \theta_i = \nu_i + \sum_{j=1}^N a_{ij} \sin(\theta_j  - \theta_i),\quad \theta\in \bbr.
\end{aligned}
\end{align}
Reduction basically follows from the property \eqref{B-4} and the projection operator, since the models \eqref{A-0-0} and \eqref{B-7-1} share the gradient flow structure. Moreover,  system \eqref{A-0-0} satisfies the following splitting property for a  homogeneous ensemble $\Xi_i \equiv \Xi$ for $i=1,\cdots,N$. Then, the rotated variable $Y_i:=   S_i e^{-t\Xi}$ satisfies \eqref{A-0-0} with $\Xi \equiv O$. For a proof on  the reduction and splitting property, we refer to Proposition 1 of \cite{Ma2}. Thus, when we consider a homogeneous ensemble, we set $\Xi \equiv O$ without loss of generality.

\begin{remark}
If we perform left-multiplication by $S_i^\top $ in \eqref{A-0-0}, then we obtain the following reduced dynamics:
\begin{equation}  \label{B-11}
S_i^\top  \dot S_i = \Xi_i + \frac{\kp}{2N} \sum_{k=1}^N \Big( S_i^\top  S_k - S_k^\top  S_i\Big),
\end{equation}
which means that \eqref{B-11} can be uniquely derived from \eqref{A-0-0}. In other words, if $\{S_i\}_{i=1}^N$ is a solution to \eqref{A-0-0}, then it also becomes a solution to \eqref{B-11}. However, the converse statement might not hold, since there does not exist an inverse matrix of $S_i^\top $. Hence, a solution set for \eqref{A-0-0} is a proper subset of that for \eqref{B-11}.
%
\end{remark}
Below, we briefly recall previous results in \cite{Ma1,Ma2} on the first-order model on the Stiefel manifolds. As far as the authors know, there has been only two literatures which concern  with the consensus model \eqref{A-1} on the Stiefel manifolds. A graph $\mathcal G$ is a pair $(\textup{V},\textup{E})$ where $\textup{V}=\{1,\cdots,N\}$ and $E$ is a subset of $\textup{V}$ where each subset is of cardinality two. Suppose that $\mathcal S =(S_1,\cdots,S_N)$ is the state of particles interacting through the graph $\mathcal G$.    
\begin{theorem} \cite{Ma1,Ma2} \label{T2.1}
Suppose the pair $(p,n)$ satisfy
\begin{equation*}
p\leq \frac23n-1,
\end{equation*}
and  $\mathcal G$ is connected. Let $\mathcal S$ be a global solution to \eqref{A-0}. Then, the consensus manifold $\mathcal C$ defined by
\begin{equation*}
\mathcal C:= \{ (S_i)_{i=1}^N \in \emph{St}(p,n)^N : S_i = S_j ,\quad \forall \{ i,j\} \in \textup E\}
\end{equation*}
is almost globally asymptotically stable.  
\end{theorem}

\subsection{A second-order extension} \label{sec:2.3}
In this subsection, we introduce a second-order extension \eqref{A-2} of \eqref{A-1} by adding suitable inertia like terms:
\begin{align} \label{B-13}
\begin{aligned}
  m  \ddot S_i &= - mS_i  \dot S_i^\top  \dot S_i - \gamma  \dot S_i + S_i \Xi_i +\frac{m}{\gamma} ( 2\dot S_i \Xi_i - S_i \Xi_i S_i^\top  \dot S_i + S_i \dot S_i^\top  S_i \Xi_i) \\
& \hspace{0.5cm}+ \frac{\kp}{N}\sum_{k=1}^N \left(  S_k - \frac12 ( S_iS_i^\top S_k + S_i S_k^\top S_i)\right).\\
\end{aligned}
\end{align}
Note that in a formal zero inertia limit $m \to 0$, system \eqref{B-13} reduces to the first-order model \eqref{A-0-0}. In order to show that the proposed model \eqref{B-13} is a suitable extension on the Stiefel manifold, we need to make sure that the governing manifold $\St$ is still positively invariant along \eqref{B-13}. 

\begin{lemma}[Positive invariance of the Stiefel manifold] \label{L2.3}
Suppose the initial data $(\mathcal S^\textup{in}, \dot{\mathcal S}^\textup{in})$ satisfy 
\begin{equation} \label{B-14}
S_i^\textup{in} \in \St,\quad \dot S_i^{t,\textup{in}} S_i^\textup{in} + S_i^{t,\textup{in}} \dot S_i^\textup{in} = O,\quad i=1,\cdots,N,
\end{equation}
and let $\mathcal S$ be a global solution to \eqref{B-13}. Then, we have
\begin{equation*}
S_i(t) \in \St,\quad i=1,\cdots,N,\quad t\geq0.
\end{equation*}

\end{lemma}

\begin{proof}
We take the  left multiplication of $S_i^\top $ to \eqref{B-13} to obtain  
\begin{align}\label{D-1}
\begin{aligned}
 m S_i^\top  \ddot S_i &= - m  S_i^\top S_i \dot S_i^\top  \dot S_i - \gamma S_i^\top  \dot S_i + S_i^\top S_i   \Xi_i +\frac{m}{\gamma} \left( 2S_i^\top  \dot S_i \Xi_i -S_i^\top S_i    \Xi_i S_i^\top \dot S_i + S_i^\top S_i  \dot S_i^\top  S_i \Xi_i  \right) 
\\
& \hspace{0.5cm}+ \frac{\kp}{2N}\sum_{k=1}^N \bigg[ S_i^\top S_k - \frac12 (S_i^\top S_i S_i^\top S_k + S_i^\top S_i S_k^\top S_i) \bigg].
\end{aligned}
\end{align}
We transpose \eqref{D-1} and use the skew-symmetry of $\Xi_i$ to find 
\begin{align}\label{D-2}
\begin{aligned}
 m \ddot S_i^\top S_i &= - m \dot S_i^\top  \dot S_i S_i^\top S_i - \gamma \dot S_i^\top  S_i - \Xi_iS_i^\top S_i  +\frac{m}{\gamma} \left(  -2\Xi_i \dot S_i^\top  S_i + \dot S_i^\top  S_i\Xi_iS_i^\top S_i  - \Xi_i S_i^\top  \dot S_i S_i^\top S_i \right) \\
&\hspace{0.5cm}+ \frac{\kp}{2N}\sum_{k=1}^N \bigg[ S_k^\top S_i - \frac12 (S_k^\top  S_i S_i^\top S_i + S_i^\top S_k S_i^\top S_i) \bigg].
\end{aligned}
\end{align}
We add \eqref{D-1} and \eqref{D-2} to get
\begin{align} \label{D-2-1}
\begin{aligned}
m( S_i^\top  \ddot S_i + \ddot S_i^\top  S_i) &= - m ( S_i^\top S_i \dot S_i^\top \dot S_i + \dot S_i^\top \dot S_i S_i^\top S_i) - \gamma (S_i^\top  \dot S_i + \dot S_i^\top S_i) + [S_i^\top S_i, \Xi_i] \\
&\hspace{0.5cm} + \frac m \gamma \mathcal J_1 +\frac{\kp}{2N}\sum_{k=1}^N  \mathcal J_{2k},
\end{aligned}
\end{align}
where $\mathcal J_1$ and $\mathcal J_{2k}$ are defined as
\begin{align*} \label{D-2-2}
\begin{aligned}
&\mathcal J_1 := 2S_i^\top  \dot S_i \Xi_i -S_i^\top S_i    \Xi_i S_i^\top \dot S_i + S_i^\top S_i  \dot S_i^\top  S_i \Xi_i -2\Xi_i \dot S_i^\top  S_i + \dot S_i^\top  S_i\Xi_iS_i^\top S_i  - \Xi_i S_i^\top  \dot S_i S_i^\top S_i, \\
&\mathcal J_{2k} := S_i^\top S_k + S_k^\top S_i  - \frac12 (S_i^\top S_i S_i^\top S_k + S_i^\top S_i S_k^\top S_i + S_k^\top  S_i S_i^\top S_i + S_i^\top S_k S_i^\top S_i).
\end{aligned}
\end{align*}
We recall the notation:
\begin{align*}
H_i  = I_p - S_i^\top S_i ,\quad \dot{H}_i = -(\dot{S}_i^\top S_i +S_i^\top \dot{S}_i)\quad\textup{and}\quad \ddot{H}_i =-( \ddot{S}_i^\top S_i +2\dot{S}_i^\top \dot{S}_i +S_i^\top \ddot{S}_i).
\end{align*}
Then, system \eqref{D-2-1} can be rewritten in terms of $H_i$ and its derivatives:
\begin{equation} \label{D-2-3}
m \ddot H_i + \gamma \dot H_i = -m (H_i \dot S_i^\top  \dot S_i + \dot S_i^\top  \dot S_i H_i ) + [H_i,\Xi_i]  -\frac m\gamma \mathcal J_1 -\frac{\kp}{2N}\sum_{k=1}^N \mathcal J_{2k}.
\end{equation}
In what follows, we present the estimates for $\mathcal J_1$ and $\mathcal J_2$, respectively. \newline

\noindent $\bullet$ (Estimate on $\mathcal J_1$): By direct calculation, one has  
\begin{align*}
\mathcal J_1 &= 2S_i^\top  \dot S_i \Xi_i - \Xi_i S_i^\top  \dot S_i + \dot S_i^\top  S_i \Xi_i + H_i( \dot S_i^\top S_i \Xi_i - \Xi_i S_i^\top  \dot S_i)  \\
&\hspace{0.5cm}-2\Xi_i \dot S_i^\top  S_i + \dot S_i^\top  S_i \Xi_i - \Xi_i S_i^\top \dot S_i - (\dot S_i^\top  S_i \Xi_i - \Xi_i S_i^\top  \dot S_i ) H_i \\
& = 2 \Xi_i \dot H_i - 2 \dot H_i \Xi_i + H_i ( \dot S_i^\top S_i \Xi_i - \Xi_i S_i^\top  \dot S_i) -   ( \dot S_i^\top S_i \Xi_i - \Xi_i S_i^\top  \dot S_i)H_i \\
&= 2 [\Xi_i, \dot H_i] + [H_i, \dot S_i^\top  S_i \Xi_i - \Xi_i S_i^\top  \dot S_i].
\end{align*}

\noindent $\bullet$ (Estimate on $\mathcal J_{2k}$): Since the communication term of first-order and second-order models are same, it follows from  Lemma \ref{L2.1} that 
\begin{align*}
\mathcal J_{2k} = H_i(S_i^\top S_k + S_k^\top S_i) + (S_i^\top S_k + S_k^\top S_i)H_i.
\end{align*}
In \eqref{D-2-3}, we use the calculation of $\mathcal J_1$ and $\mathcal J_2$ to find the second-order (autonomous) matrix-valued equation for $H_i$:
\begin{align} \label{B-2-4}
\begin{aligned}
&m \ddot H_i + \gamma \dot H_i + m H_i (\dot S_i^\top  \dot S_i) + m (\dot S_i^\top  \dot S_i) H_i - [H_i,\Xi_i] + \frac m\gamma \Big( 2 [\Xi_i, \dot H_i] + [H_i, \dot S_i^\top  S_i \Xi_i - \Xi_i S_i^\top  \dot S_i] \Big) \\
& \hspace{1cm} + \frac{\kp}{2N}\sum_{k=1}^N \Big[ H_i(S_i^\top S_k + S_k^\top S_i) + (S_i^\top S_k + S_k^\top S_i)H_i \Big] = O.
\end{aligned}
\end{align}
One can check that $H_i=O$ becomes a solution to \eqref{B-2-4} satisfying the initial assumption \eqref{B-14}. Since a solution to the Cauchy problem \eqref{B-2-4} with the initial assumption \eqref{B-14} is unique, we conclude that 
\begin{equation*}
H_i(t) = O,\quad t>0.
\end{equation*}
This yields the desired result.  
\end{proof}
Next, we focus on the situation in which solution operators for \eqref{B-13} can be expressed as a composition of two operators.  In the proof of Lemma \ref{L2.3}, we note that positive invariance of $\St$ is also valid, when the fourth term in the right-hand side of \eqref{B-13} is absent. However, the second-order model \eqref{B-13} satisfies such property for a homogeneous ensemble   when the fourth term in the right-hand side of \eqref{B-13} is included. 

\begin{lemma} \label{L2.4}
Suppose that the initial data and frequency matrices satisfy 
\begin{equation*}
S_i^\textup{in} \in \St,\quad \dot S_i^{t,\textup{in}} S_i^\textup{in} + S_i^{t,\textup{in}} \dot S_i^\textup{in} =  O,\quad \Xi_i \equiv \Xi, \quad i=1,\cdots,N, 
\end{equation*} 
and let $\mathcal S$ be a global solution to \eqref{B-13}. Then, $Y_i:= S_i e^{-\frac{\Xi}{\gamma}t}$ satisfies
\begin{equation*}
  m  \ddot Y_i = - mY_i  \dot Y_i^\top  \dot Y_i - \gamma  \dot Y_i + \frac{\kp}{N}\sum_{k=1}^N \left(  Y_k - \frac12 ( Y_iY_i^\top Y_k + Y_i Y_k^\top Y_i)\right).
\end{equation*}  
\end{lemma}

\begin{proof}
By direct calculation, we use the ansatz for $Y_i$ to find its derivatives:
\begin{align*}
&\dot Y_i = \left( \dot S_i - \frac{S_i\Xi_i}{\gamma}\right)e^{-\frac{\Xi}{\gamma}t}, \quad \ddot Y_i = \left( \ddot S_i - \frac{2 \dot S_i \Xi_i}{\gamma} + \frac{S_i \Xi_i^2}{\gamma^2} \right) e^{-\frac{\Xi}{\gamma}t}, \\
&Y_i^\top  = e^{\frac{\Xi}{\gamma}t} S_i^\top ,\quad \dot Y_i^\top  = e^{\frac{\Xi}{\gamma}t} \left( \dot S_i^\top   + \frac{\Xi_i S_i^\top }{\gamma} \right).
\end{align*}
In what follows, we consider the three terms:
\begin{equation*}
m\ddot Y_i + mY_i \dot Y_i^\top  \dot Y_i,\quad \gamma \dot Y_i,\qquad Y_k -\frac12(Y_i Y_i^\top  Y_k + Y_iY_k^\top Y_i).
\end{equation*}
$\bullet$ (Estimate of $m\ddot Y_i + mY_i \dot Y_i^\top  \dot Y_i$): we observe
\begin{align*}
 &m\ddot Y_i + mY_i \dot Y_i^\top  \dot Y_i \\
 &\hspace{0.3cm}= m\left( \ddot S_i - \frac{2 \dot S_i \Xi_i}{\gamma} + \frac{S_i \Xi_i^2}{\gamma^2}\right) e^{-\frac{\Xi}{\gamma}t} + mS_i \left( \dot S_i^\top  \dot S_i + \frac{\Xi_i S_i^\top  \dot S_i}{\gamma} - \frac{\dot S_i^\top  S_i \Xi_i}{\gamma} - \frac{\Xi_i^2}{\gamma^2}\right) e^{-\frac{\Xi}{\gamma}t} \\
 &\hspace{0.3cm} = \left( m\ddot S_i +m S_i \dot S_i^\top  \dot S_i  - \frac{2m \dot S_i \Xi_i}{\gamma} + \frac{mS_i\Xi_i S_i^\top  \dot S_i}{\gamma} - \frac{mS_i\dot S_i^\top  S_i \Xi_i}{\gamma} \right) e^{-\frac{\Xi}{\gamma}t}.
\end{align*}
$\bullet$ (Estimate of $\gamma \dot Y_i$): By direct calculation, 
\begin{equation*}
\gamma \dot Y_i = (\gamma \dot S_i - S_i\Xi_i) e^{-\frac\Xi \gamma t}.
\end{equation*}
 $\bullet$ (Estimate of $ Y_k -\frac12(Y_i Y_i^\top  Y_k + Y_iY_k^\top Y_i$)): we use the ansatz of $Y_i$ to find
 \begin{align*}
  Y_k -\frac12(Y_i Y_i^\top  Y_k + Y_iY_k^\top Y_i) = \left(  S_k -\frac12(S_iS_i^\top  S_k + S_iS_k^\top S_i)   \right) e^{-\frac\Xi\gamma t}.
 \end{align*}
Finally, we combine the calculations above to find our desired result:
\begin{align*}
  &m  \ddot Y_i +  mY_i  \dot Y_i^\top  \dot Y_i + \gamma  \dot Y_i - \frac{\kp}{N}\sum_{k=1}^N \left(  Y_k - \frac12 ( Y_iY_i^\top Y_k + Y_i Y_k^\top Y_i)\right) \\
  &\hspace{0.5cm} = \biggl(   m \ddot S_i + mS_i \dot S_i^\top  \dot S_i + \gamma \dot S_i - S_i \Xi_i - \frac m\gamma ( 2\dot S_i \Xi - S_i\Xi_i S_i^\top  \dot S_i + S_i\dot S_i^\top  S_i\Xi_i) \\
  &\hspace{0.2cm} - \frac\kp N \sum_{k=1}^N \left( S_k - \frac12(S_iS_i^\top S_k + S_iS_k^\top S_i) \right) \biggr)e^{-\frac\Xi\gamma t} =O.
\end{align*}
\end{proof}
%
We finally close this section by introducing second-order Gr\"onwall-type inequalities.

\begin{lemma}   \cite{C-H-Y, H-K1}  \label{gronwall}
Let $y=y(t)$ be a nonnegative $\mathcal C^2$-function satisfying the following differential inequality:
\begin{equation*}
a\ddot y + b \dot y + c y \leq \veps(t), \quad t>0.
\end{equation*}
Then, the following estimates hold:
\begin{enumerate}
\item Suppose that $b^2-4ac>0$ and $\veps(t) \equiv \veps_0$ is a given positive constant. Then, we have
\begin{align*}
y(t)& \leq \frac{\veps_0}{c} + \left( y(0) + \frac dc\right) e^{-\nu_1t}  \\
&+\frac{a}{\sqrt{b^2-4ac}} \left( y'(0)+ \nu_1y(0) - \frac{2\veps_0}{b-\sqrt{b^2-4ac}} \right) \left( e^{-\nu_2t}-e^{-\nu_1t}\right).
\end{align*}
where $\nu_1$ and $\nu_2$ are given as follows: 
\begin{equation*}
\nu_1 := \frac{ b+ \sqrt{b^2-4ac}}{2a}, \quad \nu_2:= \frac{b-\sqrt{b^2-4ac}}{2a}.
\end{equation*}
Then, one has
\begin{equation*}
\limsup_{t\to\infty} y(t) \leq \frac{\veps_0}{c}.
\end{equation*}
\item Suppose that $b^2-4ac<0$. Then, we have
\begin{align*}
y(t) \leq  \frac{4a\veps_0}{b^2} + \left( y(0) - \frac{4a\veps_0}{b^2} + \left( \frac{b}{2a} y(0) + y'(0) - \frac{2\veps_0}{b}\right)t \right) e^{-\frac{b}{2a}t}.
\end{align*}
Moreover, 
\begin{equation*}
\limsup_{t\to\infty} y(t) \leq \frac{4a\veps_0}{b^2}.
\end{equation*}
\item We furthermore assume that $\veps(t)$ is a nonnegative continuously differentiable function decaying to zero as $t\to\infty$. Then, we have
\begin{equation*}
\lim_{t\to\infty} y(t) = 0.
\end{equation*}
\end{enumerate}
\end{lemma}

\begin{proof}
For the first two assertions, we refer the reader to \cite[Lemma 3.1]{C-H-Y}. On the other hand for the last assertion, we refer the reader to \cite[Lemma 4.9]{H-K1}.
\end{proof}


\section{First-order consensus model} \label{sec:3}
\setcounter{equation}{0}
In this section, we study the asymptotic behavior of the first-order model. In Section \ref{sec:3.1}, we consider the homogeneous ensemble in which all natural frequencies are the same so that the complete consensus can be achieved. In Section \ref{sec:3.2}, we are concerned with the heterogeneous ensemble where that phase-locked states can arise under some framework with a large coupling strength regime.\newline

Recall the synchronization quantity $H_{ij} \in M_{p,p}(\bbr)$:
\begin{equation*}
H_{ij}:= I_p - S_i^\top  S_j,\quad i,j=1,\cdots,N.
\end{equation*}
Then, we observe from the definition of $\St$
\begin{equation} \label{C-0}
\|S_i-S_j\|_\text{F}^2 = \text{tr} ( I_p - S_i^\top  S_j + I_p - S_j^\top S_i) =  \text{tr}(H_{ij} + H_{ji}) \leq \sqrt p   \|H_{ij} + H_{ji}\|_\tF,
\end{equation}
where we used the inequality:
\begin{equation*}
\textup{tr}(AB) \leq \|A\|_\tF\|B\|_\tF \quad \textup{for $A=I_p$ and $B= H_{ij} + H_{ji}$}.
\end{equation*}
On the other hand, we see
\begin{align} \label{C-0-1}
\begin{aligned}
\|H_{ij} + H_{ji} \|_\tF &= \| I_p - S_i^\top  S_j + I_p - S_j^\top S_i \|_\tF = \|S_i^\top S_i - S_i^\top  S_j + S_j^\top  S_j - S_j^\top  S_i \|_\tF \\
& \leq \| (S_i^\top  - S_j^\top )(S_i - S_j) \|_\tF \leq \|S_i - S_j\|_\tF^2.
\end{aligned}
\end{align}
In \eqref{C-0} and \eqref{C-0-1}, one has
\begin{equation} \label{C-0-2}
 \frac{1}{\sqrt p} \|S_i - S_j\|_\tF^2 \leq \|H_{ij} + H_{ji}\|_\tF \leq \|S_i - S_j\|_\tF^2.
\end{equation}
Thus, the following equivalence holds:
\begin{equation*}
\lim_{t\to\infty} \|S_i - S_j\|_\tF =0 \quad  \Longleftrightarrow  \quad \lim_{t\to\infty}  \|H_{ij}+ H_{ji}\|_\tF = 0. 
\end{equation*}

\subsection{A homogeneous ensemble} \label{sec:3.1} 
In this subsection, we consider  a homogeneous ensemble where all frequency matrices are same:
\begin{equation*}
\Xi_i \equiv O,\quad i=1,\cdots,N.
\end{equation*}
\begin{lemma} \label{L3.1}
Let $\mathcal S$ be a solution to \eqref{A-1} with $\Xi_i \equiv \Xi$. Then, $\|S_i - S_j\|_F^2$ satisfies 
\begin{align} \label{C-0-3}
\begin{aligned}
\frac{\d}{\d t} \|S_i - S_j\|^2& \leq -\frac{\kp}{N} \sum_{k=1}^N a_{ik} ( \|S_i - S_j\|^2 - \|S_j - S_k\|^2) - \frac\kp N \sum_{k=1}^N a_{jk}( \|S_i- S_j\|^2 - \|S_i - S_k\|^2) \\
& \hspace{0.5cm}-( 2- \|S_ i- S_j\|^2) \cdot \frac{\kp}{2N } \sum_{k=1}^N( a_{ik} \|S_ i- S_k\|^2 + a_{jk} \|S_j - S_k\|^2) .
\end{aligned}
\end{align}
Moreover, for the lower bound estimate, one has
\[  \frac{\d}{\d t}d_{ij}^2 \geq -4\kp a_M  d_{ij}^2,\quad t>0.\]
\end{lemma}

\begin{proof}
First, we derive the dynamics for $ \frac{\d}{\d t}(S_j^\top  S_i) = \dot S_j^\top  S_i + S_j^\top  \dot S_i $. For this, we observe 
\begin{align} \label{C-3}
\begin{aligned}
& S_j^\top  \dot S_i =S_j^\top  S_i \Xi_i + \frac{\kp}{N}\sum_{k=1}^N a_{ik}\bigg[ S_j^\top  S_k - \frac{1}{2} ( S_j^\top  S_i S_i^\top  S_k + S_j^\top   S_i S_k^\top  S_i) \bigg], \\
&\dot S_j^\top  S_i = -\Xi_j S_j^\top S_i  + \frac{\kp}{N}\sum_{k=1}^N a_{jk}\bigg[ S_k^\top S_i  -\frac12 ( S_k^\top  S_j S_j^\top S_i  + S_j^\top  S_k S_j^\top  S_i) \bigg].
\end{aligned}
\end{align}
Then, we calculate $\eqref{C-3}_1 + \eqref{C-3}_2$ to find the dynamics for $A_{ji}:= S_j^\top  S_i$:
\begin{align}
\begin{aligned} \label{C-4}
\frac{\d}{\d t} A_{ji} &= A_{ji} \Xi_i - \Xi_j A_{ji}  \\
&\hspace{0.5cm} + \frac\kp N \sum_{k=1}^N a_{ik} \biggl[   A_{jk} - \frac12 (A_{ji} A_{ik} +A_{ji} A_{ki}) \biggl]  + a_{jk} \biggl[ A_{ki} - \frac12(A_{kj} A_{ji}  + A_{jk} A_{ji}) \biggl].
\end{aligned}
\end{align}
Again, if we recall the notation: 
\begin{equation*}
H_{ji} = I_p - A_{ji} = I_p - S_j^\top S_i,
\end{equation*}
then \eqref{C-4} directly yields the dynamics for $H_{ji}$:
\begin{align} \label{C-6}
\begin{aligned}
 \frac{\d}{\d t} H_{ji} &=  (\Xi_j - \Xi_i) + (H_{ji} \Xi_i - \Xi_j H_{ji}) -\frac\kp N \sum_{k=1}^N (a_{ik} + a_{jk}) H_{ji}  \\
& \hspace{0.5cm}+ \frac{\kp}{2N } \sum_{k=1}^N a_{ik} ( 2H_{jk} - H_{ik} - H_{ki}) + a_{ik} (H_{ji} H_{ik} + H_{ji} H_{ki}) \\
& \hspace{0.5cm} + \frac{\kp}{2N } \sum_{k=1}^N  a_{jk}( 2_{ki} - H_{kj} - H_{jk}) + a_{jk} ( H_{kj} H_{ji} + H_{jk} H_{ji}) .
\end{aligned}
\end{align}
We interchange the index $i \leftrightarrow j$ in \eqref{C-6} to  get
\begin{align} \label{C-7}
\begin{aligned}
 \frac{\d}{\d t} H_{ij} &=  (\Xi_i - \Xi_j) + (H_{ij} \Xi_j - \Xi_i H_{ij}) -\frac\kp N \sum_{k=1}^N (a_{jk} + a_{ik}) H_{ij}  \\
& \hspace{0.5cm} + \frac{\kp}{2N } \sum_{k=1}^N a_{jk} ( 2H_{ik} - H_{jk} - H_{kj}) + a_{jk} (H_{ij} H_{jk} + H_{ij} H_{kj}) \\
& \hspace{0.5cm} + \frac{\kp}{2N } \sum_{k=1}^N  a_{ik}( 2H_{kj} - H_{ki} - H_{ik}) + a_{ik} ( H_{ki} H_{ij} + H_{ik} H_{ij}) .
\end{aligned}
\end{align}
We add \eqref{C-6} and   \eqref{C-7}  to find
\begin{align}  \label{C-7-1}
\begin{aligned}
\frac{\d}{\d t} (H_{ij} + H_{ji}) & = -\frac\kp N \sum_{k=1}^N (a_{ik} + a_{jk}) (H_{ij} + H_{ji}) +  (   H_{ji} \Xi_i - \Xi_j H_{ji} + H_{ij} \Xi_j - \Xi_i H_{ij} ) \\
& \hspace{0.5cm}   + \frac{\kp}{2N} \sum_{k=1}^N (a_{jk} - a_{ik} )( H_{ik} + H_{ki} - H_{jk} - H_{kj}) \\
& \hspace{0.5cm} + \frac{\kp}{2N} \sum_{k=1}^N a_{ik} (  H_{ji}(H_{ik} + H_{ki} ) + (H_{ik}+H_{ki})H_{ij}      ) \\
&\hspace{3cm}+a_{jk} ( (H_{jk} + H_{kj})H_{ji} + H_{ij} (H_{jk} + H_{kj})).
\end{aligned}
\end{align}
We take  the trace in \eqref{C-7-1} to find 
\begin{align} \label{C-7-2}
\begin{aligned}
\frac{\d}{\d t} \|S_i - S_j\|^2& = -\frac\kp N \sum_{k=1}^N (a_{ik} + a_{jk}) \|S_i - S_j\|^2 + \textup{tr}( (H_{ij} - H_{ji})(\Xi_i - \Xi_j) ) \\
&\hspace{0.5cm} + \frac{\kp}{2N} \sum_{k=1}^N (a_{jk} - a_{ik})( \|S_i - S_k\|^2 - \|S_j - S_k\|^2) \\
&\hspace{0.5cm} + \frac{\kp}{2N}\sum_{k=1}^N \textup{tr} \Big[( a_{ik}(H_{ik} + H_{ki} )+a_{jk}( H_{jk} + H_{kj}))(H_{ij} + H_{ji}) \Big].
\end{aligned}
\end{align}
Finally, we use $\Xi_i \equiv O$ and the inequality \eqref{C-0-2} to estimate the last term in \eqref{C-7-2}: 
\begin{align*}
\frac{\d}{\d t} \|S_i - S_j\|^2& \leq -\frac{\kp}{N} \sum_{k=1}^N a_{ik} ( \|S_i - S_j\|^2 - \|S_j - S_k\|^2) - \frac\kp N \sum_{k=1}^N a_{jk}( \|S_i- S_j\|^2 - \|S_i - S_k\|^2) \\
& \hspace{0.5cm}-( 2- \|S_ i- S_j\|^2) \cdot \frac{\kp}{2N } \sum_{k=1}^N( a_{ik} \|S_ i- S_k\|^2 + a_{jk} \|S_j - S_k\|^2) .
\end{align*}
\end{proof}
We are now ready to prove a proof of Theorem \ref{T1.1}. For this, we define the maximal diameter:
\begin{equation*} \label{C-7-3}
\mathcal D(\mathcal S(t)) := \max_{1\leq i,j \leq N } \|S_i(t) - S_j (t) \|_\tF,\quad t \geq0.
\end{equation*}

\noindent \textbf{(Proof of Theorem \ref{T1.1})}:~For each $t>0$, we choose the maximal index $(i_t,j_t)$:
\begin{equation*}
D(\mathcal S(t)) = \|S_{i_t} - S_{j_t}\|_\tF. 
\end{equation*}
Then, it follows from \eqref{C-0-3} in Lemma \ref{L3.1}  that $D(\mathcal S)$ satisfies
\begin{align*}
\frac{\d}{\d t } \mathcal D(S)^2 &\leq -(2-  \mathcal D(S)^2) \cdot \frac{\kp}{2N} \sum_{k=1}^N( a_{ik} d_{i_tk}^2 + a_{jk} d_{j_tk}^2 ) \leq -\frac{\kp a_m}{4} (2-  \mathcal D(S)^2)  \mathcal D(S)^2.
\end{align*}
Thus, the region $\{ t\geq0: \mathcal D(S(t)) <\sqrt 2\}$ is positively invariant and the desired assertion directly follows from the dynamical systems theory. \qed

%
%
\begin{remark} \label{R3.1}
In Theorem \ref{T2.1}, their initial framework leading to the complete consensus  only depends on the relation of  $(p,n)$ (in fact, $p\leq \frac{2n}{3}-1$),  and the initial data $\mathcal S^\textup{in}$ is not restricted except for a measure zero set, so that the global consensus is achieved. In contrast, for our result, any $(p,n)$ can be chosen; however, as a trade-off, smallness on the initial data would be imposed, for instance, $\mathcal  D(\mathcal S^\textup{in})<\sqrt2$.  
\end{remark}

\subsection{A heterogeneous ensemble} \label{sec:3.2} 
In this subsection, we deal with the heterogeneous ensemble, and provide  a sufficient framework leading to the phase-locked state. In order to capture the phase-locked state, we briefly present our strategy introduced in \cite{H-R} consisting of  several steps: \newline

\noindent $\bullet$~Step A: let $\mathcal S$ and $\tilde{\mathcal S}$ two solutions to \eqref{A-1}. Derive a temporal evolution of the maximal quantity (inter-diameter) denoted by $d(\mathcal S,\tilde{\mathcal S})$:
\begin{equation*}
d(S,\tilde S) :=\max_{1\leq i,j\leq N} \|S_i^\top  S_j - \tilde S_i^\top \tilde S_j\|_F = \max_{1\leq i,j \leq N}\|A_{ij} - \tilde A_{ij}\|_F,
\end{equation*}
 which measures the relative distance between two solution configurations $\mathcal S$ and $\tilde{\mathcal S}$ (see Lemma \ref{L3.2}). 

\vspace{0.3cm}

\noindent $\bullet$~Step B: in Step A, since we need a smallness on the diameter $D(\mathcal S)$, we find a positively invariant region of $\mathcal D(\mathcal S)$ (see Lemma \ref{L3.3}). 

\vspace{0.3cm}

\noindent $\bullet$~Step C: together with a temporal evolution of $d(\mathcal S,\tilde{\mathcal S})$ and a positively invariant region of $\mathcal D(\mathcal S)$, we show that $d(\mathcal S, \tilde{\mathcal S})$ converges to zero. Then, since our system is autonomous, $\tilde S_i(t) = S_i(t+T)$ also becomes a solution to the system for any $T$. By discretizing the time $t\in \bbr_+$ as $n\in \bbz_+$ and setting $T=m\in \bbz_+$, we deduce that $\{S_j^\top(n)  S_i(n)\}_{n\in \bbz_+}$ is indeed a Cauchy sequence in the radius $p$-ball of $M_{p}(\bbr)$ and hence each $S_j^\top S_i$ converges to a constant $p\times p$ matrix.

\vspace{0.3cm}

In what follows, we provide the detailed justification of Step A, Step B and Step C. For this, we assume 
\begin{equation*}
  \Lambda :=   a_m - \frac{N-1}{N}(a_M + d(\mathcal A))  \in (0, 8p a_M^2).
\end{equation*}
In other words, difference between the maximum and minimum would be small. Below, we study the temporal evolution of $d(\mathcal S,\tilde{\mathcal S})$ in the following lemma.

\begin{lemma} \label{L3.2}
Let $\mathcal S$ and $\tilde{\mathcal S}$ be two solutions to \eqref{A-1}. Then, the inter-diameter $d(\mathcal S,\tilde{\mathcal S})$ satisfies
\begin{equation*} \label{C-15}
 \frac{\d}{\d t} d(\mathcal S,\tilde{\mathcal S}) \leq -2\kp (\Lambda - a_M \sqrt p (\mathcal D(S) + \mathcal D(\tilde S)))   d(\mathcal S,\tilde{\mathcal S}), \quad t>0.
\end{equation*}
\end{lemma}
\begin{proof}
 Since the proof is   lengthy, we leave its  proof in Appendix \ref{sec:A.1}. 
\end{proof}

\vspace{0.5cm}

Next, we find a positively invariant region of $\mathcal D(\mathcal S)$. For this, we define $0<\alpha<\beta$ through the relation:
\begin{equation*} \label{C-16-1}
\textup{$\alpha$ and $\beta$ are two positive roots of $r^3 - 2r + \frac{2\sqrt p \|\Xi\|_\infty}{\kp_*a_m}=0$},
\end{equation*}
where $\kp_*$ is a positive constant defined by
\[\kp_*:= \frac{16p^2 a_M^3 \|\Xi\|_\infty }{a_m(8pa_M^2 \Lambda-\Lambda^3) }. \]

\begin{lemma} \label{L3.3}
Suppose that the coupling strength and the initial data satisfy 
\begin{equation} \label{C-16}
\kp   >\max\left\{  \frac{\sqrt{6p}}{9}\frac{\|\Xi\|_\infty}{a_m},~ \frac{16p^2 a_M^3 \|\Xi\|_\infty }{a_m(8pa_M^2 \Lambda-\Lambda^3) }\right\},\quad \mathcal D(\mathcal S^\textup{in})< \beta,
\end{equation}
and let $\mathcal S$ a  global solution to \eqref{A-1}. Then, there exists a finite entrance time $T_2>0$ such that
\begin{equation*} \label{C-16-2}
\mathcal D(\mathcal S(t))< \alpha< \frac{\Lambda}{2a_M \sqrt p},\quad t>T_2.
\end{equation*} 
\end{lemma}
\begin{proof}
Recall the estimate \eqref{C-0-3} in Lemma \ref{L3.1} and find the differential inequality of $D(\mathcal S)$ for a heterogeneous ensemble:
\begin{align*}
\frac{\d}{\d t} \mathcal D(\mathcal S) & \leq - \frac{\kp a_m}{2} \mathcal D(S) + \frac{\kp a_m}{4} \mathcal D(S)^3 + \frac{\sqrt p}{2} \|\Xi\|_\infty \\
& = \frac{\kp a_m} 4 \left( -2\mathcal D(S) + \mathcal D(S)^3 + \frac{2\sqrt p \|\Xi\|_\infty}{\kp a_m}\right) =:f(\mathcal D(\mathcal S)),
\end{align*}
where we introduce an auxiliary cubic function $f$:
\begin{equation*}
f(r) =   r^3-2r + \frac{2\sqrt p  \|\Xi\|_\infty}{\kp a_m},\quad r\geq0. 
\end{equation*}
Then, it follows from the simple calculus that 
\begin{equation*}
\textup{$f$ attains the minimum at $r_*:=\sqrt{\frac23}$ with $f(r_*) =  \frac{2\sqrt p  \|\Xi\|_\infty}{\kp a_m}- \frac23\sqrt{\frac23}$}.
\end{equation*}
Thus, since $f(r_*)<0$ is  due to $\eqref{C-16}_1$,  we see that $f$ has two positive roots $r_1$ and $r_2$ such that
\begin{equation*}
0< r_1 < \sqrt{\frac23},\quad \sqrt{\frac23}<  r_2 <\sqrt 2.   
\end{equation*}
Note that $\eqref{C-16}_1$ indeed gives $f(r_*)<0$.  In addition, we observe  
\begin{equation*}
\lim_{\kp\to\infty} r_1 = 0 \quad \lim_{t\to\infty} r_2 = \sqrt2. 
\end{equation*}
We use the dynamical systems theory to show that if $\mathcal D(\mathcal S^\textup{in}) < r_2$, then there exists a finite entrance time $T_3>0$ such that
\begin{equation*}
\mathcal D(\mathcal S(t))< r_1,\quad t>T_3.
\end{equation*}
More precisely, we split the case into two parts. First, suppose initial data satisfy 
\[ \mathcal D(S^\textup{in}) \leq r_1. \]
Then, at the time $t=t_*$ when $\mathcal D(S(t_*)) = r_1$, we have 
\begin{equation*}
\frac{\d}{\d t}\mathcal  D(S) \leq 0.
\end{equation*}
Thus, $\mathcal D(S)$ does not increase at $t=t_*$ and hence is restricted in the interval $[0,r_1]$. Second, suppose that the initial data satisfy $ r_1 <\mathcal D(S^\textup{in}) <r_2$. Again, at the instant time $t=t_*$ when $r_1 <\mathcal D(S(t_*)) < r_2$, 
\begin{equation*}
\frac{\d}{\d t} \mathcal D(S) < 0.
\end{equation*}
Thus, $\mathcal D(S)$ starts to strictly decrease. Then, by applying same argument in Proposition 3.1 in \cite{C-H-J-K}, we find such finite entrance time $T_3>0$. \newline

In order to find the desired invariant region, we have to assume 
\begin{equation}\label{C-18}
f\left( \frac{\Lambda}{2a_M\sqrt p}\right) = \left( \frac{\Lambda}{2a_M\sqrt p}\right)^3 - 2\cdot \left( \frac{\Lambda}{2a_M\sqrt p}\right) + \frac{2\sqrt p \|\Xi\|_\infty}{\kp a_m} <0,
\end{equation}
and one can check from algebraic manipulation that the relation \eqref{C-18} is equivalent to $\eqref{C-16}_1$. Then, the cubic equation 
\begin{equation*}
r^3 - 2r + \frac{2\sqrt p \|\Xi\|_\infty}{\kp_*a_m} =0
\end{equation*}
has two positive roots:  $\alpha$ and $\beta$. Hence, under the assumption \eqref{C-16}, there exists such finite entrance time $T_2>0$ so that
\begin{equation*}
\mathcal D(\mathcal S(t))<\alpha< \frac{\Lambda}{2a_M\sqrt p},\quad t>T_2.
\end{equation*} 
\end{proof}


Finally, we use Lemma \ref{L3.2} and Lemma \ref{L3.3} to present a proof of Theorem \ref{T1.2}. \newline 

%
%
\noindent \textbf{(Proof of Theorem \ref{T1.2})}:~Since we assume \eqref{C-16}, we use Lemmas \ref{L3.2} and \ref{L3.3} to show that
\begin{equation*}
\frac{\d}{\d t} d(\mathcal S, \tilde{\mathcal S})  \leq - 2\kp ( \Lambda - 2a_M\sqrt p\alpha) d(\mathcal S, \tilde{\mathcal S}),\quad \textup{a.e. $t>T_2$.}
\end{equation*}
This implies 
\begin{equation} \label{C-19}
d(\mathcal S, \tilde{\mathcal S}) (t) \leq d(\mathcal S, \tilde{\mathcal S})(T_2) e^{- 2\kp ( \Lambda - 2a_M\sqrt p\alpha) (t-T_2)},\quad \textup{a.e. $t>T_2.$}
\end{equation}
Since we are only interested in large-time behavior, without loss of generality, we would set $T_2=0$. On the other hand for any $T\geq0$, a shifted solution $\{S_i(t+T)\}$ also becomes a global solution to \eqref{A-1} with a shifted initial data $\{S_i(T)\}$, as \eqref{A-1} is autonomous. Then, \eqref{C-19} yields
\begin{equation} \label{C-20}
\|S_j^\top S_i(t+T) - S_j^\top S_i(t) \| \leq \max_{1\leq i,j\leq N} \|S_j^\top S_i(T) - S_j^\top S_i(0) \| e^{- 2\kp ( \Lambda - 2a_M\sqrt p\alpha) t },\quad t\geq0.
\end{equation}
Especially for $T=1$ and $t=n \in \bbz_+$ in \eqref{C-20}, we have
\begin{equation} \label{C-21}
\|S_j^\top S_i(n+1) - S_j^\top S_i(n) \| \leq \max_{1\leq i,j\leq N} \|S_j^\top S_i(1) - S_j^\top S_i(0) \| e^{- 2\kp ( \Lambda - 2a_M\sqrt p\alpha) n },\quad n\in \bbz_+.
\end{equation}
By induction argument, \eqref{C-21} gives for $m \in \bbz_+$, 
\begin{equation*}
\|S_j^\top S_i(n+m) - S_j^\top S_i(n) \| \leq \max_{1\leq i,j\leq N} \|S_j^\top S_i(1) - S_j^\top S_i(0) \| \frac{e^{- 2\kp ( \Lambda - 2a_M\sqrt p\alpha) n }}{1- e^{ - 2\kp ( \Lambda - 2a_M\sqrt p\alpha) }},\quad n\in \bbz_+.
\end{equation*}
Hence, the discretized sequence $\{ S_j^\top (n)S_i(n)\}_{n \bbz_+}$ is indeed Cauchy in the $p$-ball $B_p(O):=\{ X \in M_{p}(\bbr): \|X\|_\textup{F} \leq p \}$. Consequently, for each $i,j$, it converges to a constant $p\times p $ matrix $\Gamma_{ji}^\infty \in B_p(O)$. \qed

\vspace{0.2cm}

As a corollary, we show that normalized velocities synchronize. 

\begin{corollary}
Suppose that the coupling strength and initial data satisfy \eqref{C-16} and let $\mathcal S$  be a global solution to \eqref{A-1}. Then, the normalized velocities $S_i^\top  \dot S_i$ and $S_j^\top  \dot S_j$ synchronize:
\begin{equation*}
\lim_{t\to\infty} \|S_i^\top  \dot S_i - S_j^\top  \dot S_j\|_\tF = 0.
\end{equation*}

\end{corollary}

\begin{proof}
For two solutions $\mathcal S$ and $\tilde{\mathcal S}$, we observe from \eqref{B-11}
\begin{align*}
S_i^\top  \dot S_i - \tilde{S}^\top _i\dot{\tilde S}_i &=\left[ \Xi_i + \frac{\kp}{2N}\sum_{k=1}^Na_{ik}  (S_i^\top S_k - S_k^\top  S_i) \right] - \left[ \Xi_i + \frac{\kp}{2N}\sum_{k=1}^N a_{ik} (\tilde S_i^\top \tilde S_k - \tilde S_k^\top  \tilde S_i) \right] \\
&= \frac{\kp}{2N} \sum_{k=1}^N  a_{ik} \Big[ (S_i^\top S_k- \tilde S_i^\top \tilde S_k) -( S_k^\top  S_i - \tilde S_k^\top  \tilde S_i) \Big].
\end{align*}
Thus, we find
\begin{align*}
\| S_i^\top  \dot S_i - \tilde{S}^\top _i\dot{\tilde S}_i \|_\tF \leq \frac{\kp}{2N} \sum_{k=1}^N a_{ik} ( \|S_i^\top S_k- \tilde S_i^\top \tilde S_k \| + \|S_k^\top  S_i - \tilde S_k^\top  \tilde S_i\| ) \leq \kp a_M d(\mathcal S,\tilde{\mathcal S}).
\end{align*}
Finally, we use the same argument in Theorem \ref{T1.2} to show that $\{S_i^\top  \dot S_i\}$ becomes a Cauchy sequence and hence it converges to the same constant matrix. This shows the desired result. 
\end{proof}



\section{A second-order consensus model} \label{sec:4}
\setcounter{equation}{0}
In this section, we study the emergent dynamics of the second-order model \eqref{A-2}. Before we present the estimates, energy functionals are defined and their temporal evolutions are derived as a crucial step. After that, the homogeneous ensemble are considered in Section \ref{sec:4.1} and we provide a sufficient framework leading to the complete consensus. On the other hand, we deal with the heterogeneous ensemble in Section \ref{sec:4.2} and present a sufficient framework for the practical consensus in the large coupling strength and small inertia regime. \newline

As in \cite{H-K2}, we define a total energy functional associated to \eqref{A-2}:
\begin{align*}
\cE =: \underbrace{\frac{m}{N} \sum_{i=1}^N \|\dot{S}_i \|^2_\tF}_{=: \cK} + \underbrace{\frac{\kappa}{2N^2} \sum_{i,j=1}^N a_{ij} \left\|S_i-S_j\right\|_\tF^2}_{=: \cL},
\end{align*}
where $\mathcal K$ represents a (total) kinetic energy and $\mathcal L$ describes the interaction energy between the states of particles $\{S_i\}$ measuring the degree of consensus. 

\vspace{0.3cm}

Next, we study the temporal-evolution of the total energy $\cE$. First, we differentiate $\| \dot{S}_i \|_\tF^2$ with respect to $t$ to find 
\begin{align*}
m\frac{\d}{\d t}\|\dot{S}_i \|^2_\tF = m\frac{\d}{\d t}\textup{tr} ( \dot{S}_i^\top \dot{S}_i  ) = m \textup{tr}( \ddot{S}_i^\top \dot{S}_i +\dot{S}_i^\top \ddot{S}_i ) = m\textup{tr} ( \dot{S}_i^\top \ddot{S}_i + ( \dot{S}_i^\top \ddot{S}_i  )^\top   ).
\end{align*}
First, we observe the term $m\dot{S}_i^\top \ddot{S}_i$:
\begin{align} \label{D-10}
\begin{aligned}
\displaystyle m\dot{S}_i^\top \ddot{S}_i &= - m\dot{S}_i^\top S_i  \dot S_i^\top  \dot S_i - \gamma \dot{S}_i^\top  \dot S_i + \dot{S}_i^\top S_i \Xi_i + \frac{m}{\gamma} (2 \dot S_i^\top  \dot S_i \Xi_i - \dot S_i^\top  S_i \Xi_i S_i^\top  \dot S_i + \dot S_i^\top  S_i \dot S_i^\top  S_i \Xi_i) \\
& \hspace{0.5cm} +\frac{\kappa}{N}\sum_{k=1}^N a_{ik} \left( \dot{S}_i^\top  S_k - \frac12 ( \dot{S}_i^\top S_iS_i^\top S_k + \dot{S}_i^\top S_i S_k^\top S_i)\right) \\
&=: m\mathcal{I}_1-\gamma\mathcal{I}_2+\mathcal{I}_3 +\frac{m}{\gamma}\mathcal{I}_4+\frac{\kappa}{N}\sum_{k=1}^N\mathcal{I}_5.
\end{aligned}
\end{align}
In next lemma, we provide estimates for $\mathcal I_k,~k=1,\cdots,5$.
\begin{lemma}
Let $\cI_k,~ k=1,\cdots,5$ be the terms defined in \eqref{D-10}. Then, one has 
\begin{align*}
\begin{aligned} 
& \textup{tr}\left( \cI_1 +\cI_1^\top  \right) = 0,\quad \textup{tr}\left( \cI_2 +\cI_2^\top  \right) = 2\|\dot{S}_i\|_\tF^2,\quad \textup{tr}\left( \cI_3 +\cI_3^\top  \right) = \textup{tr} \left( \dot{S}_i^\top S_i\Xi_i -\Xi_i S_i^\top \dot{S}_i \right) \\
& \textup{tr}\left( \cI_4 +\cI_4^\top  \right) = 0,\quad \textup{tr}\left( \cI_5 +\cI_5^\top  \right) = a_{ik}  \textup{tr} \left( S_i^\top S_k +S_k^\top S_i \right).
\end{aligned}
\end{align*}
\end{lemma}

\begin{proof}
Below, we present estimates  of $\textup{tr}(\cI_k+{\mathcal I}_k^\top ),~ k=1,\cdots,5$, respectively. \newline 

\noindent $\bullet$ (Estimate of $\textup{tr}( \cI_1 +\cI_1^\top   )$): we recall the identity $\dot{S}_i^\top S_i +S_i^\top \dot{S}_i = O_p$ to see
\begin{align*}
\textup{tr}\left( \cI_1+\cI_1^\top  \right) = \textup{tr}( \dot{S}_i^\top S_i  \dot S_i^\top  \dot S_i +\dot{S}_i^\top \dot{S}_iS_i^\top \dot{S}_i ) = \textup{tr}\left[ \dot S_i^\top  \dot S_i \left( \dot{S}_i^\top S_i +S_i^\top \dot{S}_i \right)\right] = 0.
\end{align*}
\noindent $\bullet$  (Estimate  of $\textup{tr}\left( \cI_2 +\cI_2^\top  \right)$): by direct observation, 
\begin{align*}
\textup{tr}\left(\cI_2+\cI_2^\top  \right) = 2\textup{tr} \left( \cI_2 \right) = 2\|\dot{S}_i \|^2_\tF.
\end{align*}
\noindent $\bullet$  (Estimate of $\textup{tr}\left( \cI_3 +\cI_3^\top  \right)$): we use the property of the trace $\textup{tr}(AB) = \textup{tr}(BA)$ to find 
\begin{align*}
\textup{tr}( \cI_3 +\cI_3^\top  ) = \textup{tr} \left( \dot{S}_i^\top S_i\Xi_i -\Xi_i S_i^\top \dot{S}_i \right).
\end{align*}
\noindent $\bullet$  (Estimate of $\textup{tr}\left( \cI_4 +\cI_4^\top  \right)$): we first see
\begin{align} \label{D-11}
\begin{aligned}
\textup{tr} ( \mathcal I_4 ) & = \textup{tr}\Big(    (  2\dot S_i^\top  \dot S_i - S_i^\top  \dot S_i \dot S_i^\top  S_i + \dot S_i^\top  S_i \dot S_i^\top  S_i    )\Xi_i        \Big) \\
&= \textup{tr}\Big(    (  2\dot S_i^\top  \dot S_i - S_i^\top  \dot S_i \dot S_i^\top  S_i - \dot S_i^\top  S_i S_i^\top  \dot S_i     )\Xi_i        \Big).
\end{aligned}
\end{align}
Since the matrix $2\dot S_i^\top  \dot S_i - S_i^\top  \dot S_i \dot S_i^\top  S_i - \dot S_i^\top  S_i S_i^\top  \dot S_i $ in \eqref{D-11} is symmetric and $\Xi_i$ is skew-symmetric,   we obtain
\begin{align*}
\textup{tr}(\mathcal I_4 + {\mathcal I}_4^\top ) =0. 
\end{align*}
\noindent $\bullet$  (Estimate of $\textup{tr}\left( \cI_5 +\cI_5^\top  \right)$): we again use the identity $\dot{S}_i^\top S_i +S_i^\top \dot{S}_i = 0$ to see
\begin{align*}
\begin{aligned}
& \textup{tr}( \cI_5 +\cI_5^\top  ) \\
 &\hspace{0.5cm}=  a_{ik} \textup{tr}( \dot{S}_i^\top S_k +S_k^\top \dot{S}_i )-\frac{a_{ik}}{2}\textup{tr}( \dot{S}_i^\top S_iS_i^\top S_k + \dot{S}_i^\top S_i S_k^\top S_i +S_k^\top S_iS_i^\top \dot{S}_i +S_i^\top S_kS_i^\top \dot{S}_i ) \\
&\hspace{0.5cm}= a_{ik} \textup{tr}( \dot{S}_i^\top S_k +S_k^\top \dot{S}_i )-\frac{ a_{ik}}{2}\textup{tr}\big[ \dot{S}_i^\top S_i ( S_i^\top S_k +S_k^\top S_i  ) +\left( S_k^\top S_i +S_i^\top S_k \right) S_i^\top \dot{S}_i \big] \\
&\hspace{0.5cm}=   a_{ik}\textup{tr}( \dot{S}_i^\top S_k +S_k^\top \dot{S}_i )-\frac{ a_{ik}}{2}\textup{tr}\big[ ( \dot{S}_i^\top S_i +S_i^\top \dot{S}_i )( S_i^\top S_k +S_k^\top S_i ) \big] =  a_{ik} \textup{tr}( \dot{S}_i^\top S_k +S_k^\top \dot{S}_i).
\end{aligned}
\end{align*}
\end{proof}
Next, we provide an estimate on the time-rate of change of the total energy. 
\begin{proposition} \label{P4.1}
Let $\cS$ be a  global solution to \eqref{A-2}. Then, the total energy functional $\cE$ satisfies
\begin{align}\label{D-24}
\frac{\d\cE}{\d t} = -\frac{2\gamma}{N}\sum_{i=1}^N\| \dot{S}_i \|_\tF^2 +\frac{1}{N}\sum_{i=1}^N\textup{tr} ( \dot{S}_i^\top S_i\Xi_i -\Xi_i S_i^\top \dot{S}_i ).
\end{align}
\end{proposition}

\begin{proof}
We collect all the estimates of $\textup{tr}( \cI_i +\cI_i^\top  )$ for $i=1,\cdots,5$ to find
\begin{equation}\label{D-21}
m\frac{\d}{\d t}\|\dot{S}_i  \|^2_\tF = -2\gamma\| \dot{S}_i \|_\tF^2+\textup{tr} ( \dot{S}_i^\top S_i\Xi_i -\Xi_i S_i^\top \dot{S}_i ) +\frac{\kappa}{N}\sum_{k=1}^N  a_{ik} \textup{tr}(\dot{S}_i^\top S_k +S_k^\top \dot{S}_i ), 
\end{equation}
and   sum \eqref{D-21} over all $i$ and divide the resulting relation by $N$ to obtain
\begin{align}\label{D-22}
\begin{aligned}
\frac{m}{N}&\frac{\d}{\d t}\sum_{i=1}^N  \|\dot{S}_i \|^2_\tF = -2\frac{\gamma}{N}\sum_{i=1}^N\| \dot{S}_i \|_\tF^2 +\frac{1}{N}\sum_{i=1}^N\textup{tr} ( \dot{S}_i^\top S_i\Xi_i -\Xi_i S_i^\top \dot{S}_i ) +\kappa\textup{tr} (\dot{S}_c^\top S_c +S_c^\top \dot{S}_c ) \\
&= -2\frac{\gamma}{N}\sum_{i=1}^N\| \dot{S}_i \|_\tF^2 +\frac{1}{N}\sum_{i=1}^N\textup{tr} ( \dot{S}_i^\top S_i\Xi_i -\Xi_i S_i^\top \dot{S}_i ) +\frac{\kappa}{N^2} \sum_{i,k=1}^N a_{ik} \textup{tr} ( \dot S_i ^\top  S_k + S_k^\top  \dot S_i ) .
\end{aligned}
\end{align}
On the other hand, we observe
\begin{align*}
&\frac{\kappa}{N^2} \sum_{i,k=1}^N a_{ik} \textup{tr} ( \dot S_i ^\top  S_k + S_k^\top  \dot S_i ) = \frac{\kp}{2N^2} \sum_{i,k=1}^N a_{ik} \textup{tr} ( \dot S_i^\top  S_k + S_k^\top  \dot S_i + \dot S_k^\top  S_i + S_i^\top  \dot S_k) \\
&=\frac{\kp}{2N^2}\frac{\d}{\d t} \left( \sum_{i,k=1}^N a_{ik} \textup{tr} ( S_i^\top  S_k + S_k^\top  S_i)\right)  = -\frac{\kp}{2N^2}\frac{\d}{\d t}\left( \sum_{k=1}^N a_{ik} \|S_i - S_k\|^2\right).
\end{align*}
Finally in \eqref{D-22}, we obtain the desired  dissipation energy estimate.
\end{proof}

\begin{remark} \label{rmk:D-1}
It follows from \eqref{D-21} that
\begin{align*}
m\frac{\d}{\d t}\| \dot{S}_i \|_\tF^2 \leq -2\gamma\| \dot{S}_i \|_\tF^2 +2\|\dot{S}_i\|_\tF\cdot\left\|S_i\Xi_i\right\|_\tF +2\ka a_M \sqrt{p}\|\dot{S}_i\|_\tF, \quad t>0, 
\end{align*}
or equivalently
\begin{align*}
m\frac{\d}{\d t}\|\dot{S}_i\|_\tF \leq -\gamma\|\dot{S}_i\|_\tF +\left\|\Xi_i\right\|_\tF +\ka a_M \sqrt{p},\quad t>0.
\end{align*}
We set the maximal value of $\|\Xi_i\|_\tF$:
\begin{equation*}
\|\Xi\|_\infty:= \max_{1\leq i \leq N } \|\Xi_i\|_\tF.
\end{equation*}
Then, Gr\"onwall's lemma yields a uniform boundedness and uniform continuity of $\|\dot{S}_i\|_\tF$:
\begin{align}\label{D-24.5}
\sup_{0\leq t<\infty} \| \dot{S}_i(t) \|_\tF \leq \max \left\{ \max_{1 \leq i \leq N} \| \dot{S}_i^\textup{in} \|_\tF, \frac{1}{\gamma} \big( \|\Xi\|_\infty +\kappa a_M \sqrt{p}\big) \right\}.
\end{align}
\end{remark}

\vspace{0.5cm}

\noindent In what follows, we impose the condition on the network topology:
\begin{equation*}
\frac1N \sum_{k=1}^N a_{ik} \equiv \xi >0,\quad i=1,\cdots,N.
\end{equation*}
In other words, every columns (or rows) have a common average $\xi$. We denote several quantities:
\begin{align*}
\mathcal G(t) := \frac{1}{N^2} \sum_{i,j=1}^N \left\| S_i(t) -S_j(t) \right\|_\tF^2,\quad \mathcal D( \dot{\cS}(t)) := \max_{1 \leq i \leq N} \| \dot{S}_i(t) \|_\tF, \quad \|\Xi_i\|_\infty := \max_{1 \leq i \leq N} \left\| \Xi_i \right\|_\tF.
\end{align*}
Below, we derive a second-order differential inequality for $\mathcal G$. 
\begin{lemma} \label{L4.2} 
Let $\cS$ be a solution to \eqref{A-2}. Then, $\mathcal G$ satisfies
\begin{align}\label{D-25}
\begin{aligned}
m \ddot {\mathcal G} + \gamma \dot{\mathcal G} + 2\kp \xi \mathcal G  \leq 16m \mathcal  D(\dot{\cS})^2 +8\|\Xi\|_\infty    +\frac{16m\sqrt{p}\|\Xi\|_\infty}{\gamma}\mathcal D(\dot{\cS}).
\end{aligned}
\end{align}
\end{lemma}
\begin{proof}
Since the proof is lengthy, we postpone the justification of \eqref{D-25} in Appendix \ref{sec:B}. 
\end{proof}

Below, we quote Barbalat's lemma without proofs.
\begin{lemma}  \label{L4.2}
\emph{\cite{Ba}}
\emph{(i)}~Suppose that  a real-valued function $f: [0, \infty) \to \bbr$ is uniformly continuous and satisfies
\[ \lim_{t \to \infty} \int_0^t  f(s) \,\d s \quad \textup{exists}. \]
Then, $f$ tends to zero as $t \to \infty$:
\[ \lim_{t \to \infty} f(t) = 0. \]
\emph{(ii)}~Suppose that a real-valued function $f: [0, \infty) \to \bbr$ is continuously differentiable, and $\lim_{t \to \infty} f(t) = \alpha \in \bbr$. If $f^{\prime}$ is uniformly continuous, then 
\[ \lim_{t \to \infty} f^{\prime}(t)  = 0. \]
\end{lemma}

\subsection{A homogeneous ensemble} \label{sec:4.1} 
In this part, we consider a homogeneous ensemble whose frequency matrices are the same, i.e., $\Xi_i \equiv \Xi$.  Due to the  property in Lemma \ref{L2.4} for a homogeneous ensemble, without loss of generality, we may set 
\[ \Xi_i \equiv O, \quad i = 1, \cdots, N. \]
We are now ready to  provide  a proof of Theorem \ref{T1.3}. \newline

\noindent \textbf{(Proof of Theorem \ref{T1.3})} 
For the zero convergence of $\|\dot S_i\|_\tF$, we use in \eqref{D-24} the condition $\Xi_i \equiv O$  to get
\begin{align*}
\frac{\d\cE}{\d t} = -\frac{2\gamma}{N} \sum_{i=1}^N \| \dot{S}_i \|_\tF^2,\quad t>0.
\end{align*}
We integrate the above equation on the interval $[0,t]$ to obtain
\begin{align*}
\frac{2\gamma}{N} \int_{0}^t  \sum_{i=1}^N \| \dot{S}_i (s) \|_\tF^2\, \d s = \cE(0)-\cE(t) \leq \cE(0).
\end{align*}
Hence, one has
\begin{align*}
\int_{0}^\infty \| \dot{S}_i (t) \|_\tF^2\, \d t < \infty.
\end{align*}
It follows from Remark \ref{rmk:D-1} that $\left| \frac{\d}{\d t} \| \dot{S}_i \|_\tF \right|$ is uniformly bounded, which implies the uniform continuity of $\| \dot{S}_i \|_\tF$. Hence, we can apply Barbalat's lemma \cite{Ba} to derive the desired zero convergence of $\|\dot S_i\|_\tF$. On the other hand for the zero convergence of $\mathcal G$,  we recall the differential inequality \eqref{D-25}. Then, since we have $\|\Xi\|_\infty = 0$ and $\mathcal D(\dot S)$ converges to zero, we apply Lemma \ref{gronwall} to obtain the desired result. \qed

\subsection{A heterogeneous ensemble} \label{sec:4.2}
In this part, we study a heterogeneous ensemble and present a proof of Theorem \ref{T1.4}. For this, we impose a small inertia assumption such that there exist $m_0>0$ which does not depend on $\kp$ and $\eta>0$ satisfying
\begin{equation*}
m = \frac{m_0}{\kp^{1+\eta}}.
\end{equation*}
Since $\kp$ will tend to infinity, the relation above is referred as a small inertial assumption. The reason why we need such assumption is clarified throughout the proof. \newline

\noindent \textbf{(Proof of Theorem \ref{T1.4})}: since we assume $\eqref{Y-1}_1$, it follows from \eqref{D-24.5} that 
\begin{equation*}
\mathcal D(\dot S(t)) <  \frac1\gamma ( \|\Xi\|_\infty + \kp a_M \sqrt p ).
\end{equation*}
Then, the right-hand side of \eqref{D-25} becomes
\begin{align*}
& 16mD(\dot{\cS})^2 +8\|\Xi\|_\infty    +\frac{16m\sqrt{p}\|\Xi\|_\infty}{\gamma}\mathcal D(\dot{\cS}) \\
& \hspace{0.5cm} \leq \frac{16m}{\gamma^2} (\|\Xi\|_\infty + \kp a_M \sqrt p )^2 + \frac{16m \sqrt p \|\Xi\|_\infty} {\gamma^2} (\|\Xi\|_\infty + \kp a_M \sqrt p ) + 8 \|\Xi\|_\infty \\
& \hspace{0.5cm}= \frac{16m}{\gamma^2}(\|\Xi\|_\infty + \kp a_M \sqrt p )(\|\Xi\|_\infty + \kp a_M \sqrt p  + \sqrt p \|\Xi\|_\infty) + 8 \|\Xi\|_\infty \\
& \hspace{0.5cm} = \mathcal O(m\kp^2) + 8\|\Xi\|_\infty = \mathcal O (m\kp^2)  +\mathcal O(1).
\end{align*}
In order to apply second-order Gr\"onwall's inequality in Lemma \ref{gronwall} is used,  we consider two cases:
\begin{equation} \label{cond}
\textup{(i)}~~\gamma^2 - 8\xi m\kp >0. \qquad \textup{(ii)}~~\gamma^2- 8\xi m \kp <0.
\end{equation}
For the first case, we have
\begin{equation} \label{Y-5}
\limsup_{t\to\infty} \mathcal G(t) \leq \frac{1}{2\xi \kp} ( \mathcal O(m\kp^2) + \mathcal O(1)) = \mathcal O(m\kp) + \mathcal O(\kp^{-1}) =\mathcal O(\kp^{-\eta}) + \mathcal O(\kp^{-1}),
\end{equation}
where we used the relation $\eqref{Y-1}_2$ in the last equality. On the other hand for the second case, we have
\begin{equation} \label{Y-6}
\limsup_{t\to\infty} \mathcal G(t) \leq \frac{4m}{\gamma^2} ( \mathcal O(m\kp^2) + \mathcal O(1))  = \mathcal O(m^2\kp^2) + \mathcal O(m) = \mathcal O(\kp^{-2\eta}) + \mathcal O(\kp^{-(1+\eta)}),
\end{equation}
where $\eqref{Y-1}_2$ is used. By combining \eqref{Y-5} and \eqref{Y-6} and letting $\kp$ to infinity, we derive the desired estimate. 
\qed

%
%

\begin{remark}
Then, two conditions in \eqref{cond} can be rewritten in terms of $m_0$ and $\eta$:
\begin{equation*}
\textup{(i)}~~\kp > \left( \frac{8\xi m_0}{\gamma^2}\right)^\frac1\eta. \qquad \textup{(ii)}~~\gamma < \left(\frac{8\xi m_0}{\kp^\eta}\right)^\frac12.
\end{equation*}
In (i), we only require a large coupling strength $\kp$ as $\kp\to\infty$ and there is no assumption on $\gamma$. However for (ii), $\gamma$ depends on $\kp$. In fact, $\gamma$ also tends to zero when $\kp\to\infty$ as $m$ does. To thi

\end{remark}

\section{Conclusion}  \label{sec:5}
\setcounter{equation}{0}
In this paper, we have studied emergent dynamics of  the first-order and second-order consensus models with both homogeneous and heterogeneous ensembles in which all states of the particles are restricted to the Stiefel manifold $\St$  consisting of all orthonormal $k$-frames in $\bbr^n$. As mentioned in Section \ref{sec:1}, due to the structural property of Stiefel manifolds, it has been widely used for  optimization problems especially in the engineering community. Our first-order model  was first suggested in \cite{Ma1,Ma2}   based on the gradient flow of the total distance between all particles. For a homogeneous ensemble in which frequency matrices are the same, we show that the relaxation process towards the complete consensus state is always achieved  with an exponential rate. In addition, we provide a subset of initial data leading to the complete consensus independent of the choice $(p,n)$, whereas the results in \cite{Ma1,Ma2} show that the complete consensus is globally stable if $(p,n)$ satisfies the specific relation and hence, the global consensus is achieved. On the other hand, for the heterogeneous ensemble, the phase-locked state can emerge under the large coupling strength regime, when the smallness assumption on the initial data is imposed. We also study the asymptotic emergent dynamics of a second-order consensus model which naturally extends the first-order model by incorporating the inertial force. For a homogeneous ensemble , we present a sufficient framework leading to the complete consensus state based on a second-order Gr\"onwall's inequality and energy estimate. In contrast, for  a heterogeneous ensemble, a sufficient framework for the practical consensus state is provided under the large coupling and small inertia regime. In fact, the phase-locked state for the second-order model is not considered in this work. Thus, we leave this problem for a  future work.


\appendix

\section{Proof of Lemma \ref{L3.2}}  \label{sec:A.1}
\setcounter{equation}{0} 
In this appendix, we provide a  proof of Lemma \ref{L3.2} in which the differential inequality for $d(\mathcal S, \tilde{\mathcal S})$ is derived. \newline 

 Recall the relation \eqref{C-4} in Lemma \ref{L3.1}
 \begin{align*}
\begin{aligned}  
\frac{\d}{\d t} A_{ji} &= A_{ji} \Xi_i - \Xi_j A_{ji}  \\
&\hspace{0.5cm} + \frac\kp N \sum_{k=1}^N a_{ik} \biggl[   A_{jk} - \frac12 (A_{ji} A_{ik} +A_{ji} A_{ki}) \biggl]  + a_{jk} \biggl[ A_{ki} - \frac12(A_{kj} A_{ji}  + A_{jk} A_{ji}) \biggl] \\
& =  A_{ji} \Xi_i - \Xi_j A_{ji}  + \frac{\kp}{2N} \sum_{k=1}^N a_{ik} ( A_{jk} - A_{kj}) + a_{jk}(A_{ki} - A_{ik})  \\
&\hspace{0.5cm} + \frac{\kp}{2N} \sum_{k=1}^N \Big(a_{ik} ( A_{jk} + A_{kj}) + a_{jk}(A_{ki} + A_{ik})\\
&\hspace{3cm} - a_{ik} A_{ji}(A_{ik} + A_{ki}) - a_{jk}(A_{jk} +A _{kj})A_{ji}\Big). 
\end{aligned}
\end{align*}
Similarly, we find
\begin{align*}
\begin{aligned}
\frac{\d}{\d t} \tilde A_{ji} &= \tilde A_{ji} \Xi_i - \Xi_j \tilde A_{ji} +  \frac{\kp}{2N} \sum_{k=1}^N a_{ik} (\tilde A_{jk} - \tilde A_{kj}) + a_{jk}(\tilde A_{ki} - \tilde A_{ik})  \\
& \hspace{0.5cm}+ \frac{\kp}{2N} \sum_{k=1}^N \Big(a_{ik} (\tilde  A_{jk} +\tilde  A_{kj}) + a_{jk}(\tilde A_{ki} + \tilde A_{ik}) \\
&\hspace{3cm} - a_{ik} \tilde A_{ji}(\tilde A_{ik} + \tilde A_{ki}) - a_{jk}(\tilde A_{jk} +\tilde A _{kj})\tilde A_{ji} \Big).
\end{aligned}
\end{align*}
Below, we consider the dynamics of $A_{ji} - \tilde A_{ji}$:
\begin{align} \label{C-12-2}
\begin{aligned}
\frac{\d}{\d t} (A_{ji} - \tilde A_{ji}) &= (A_{ji} - \tilde A_{ji})\Xi_i - \Xi_j(A_{ji} - \tilde A_{ji} ) + \frac{\kp}{2N} \sum_{k=1}^N \mathcal J_{3k} \\
&\hspace{0.5cm}+\frac{\kp}{2N} \sum_{k=1}^N \Big(a_{ik}(( A_{jk} - \tilde A_{jk}) - (A_{kj} - \tilde A_{kj}))  \\
&\hspace{3cm}+ a_{jk}(( A_{ki} - \tilde  A_{ki}) - (A_{ik} - \tilde A_{ik})) \Big),
\end{aligned}
\end{align}
where $\mathcal J_3$ is defined as 
\begin{align} \label{C-12-3}
\begin{aligned}
\mathcal J_{3k} &: = \big( a_{ik} ( A_{jk} - \tilde A_{jk}) - a_{jk}(A_{jk} A_{ji} - \tilde A_{jk} \tilde A_{ji}) \big)  \\
&\hspace{0.5cm}+ \big( a_{ik} ( A_{kj} - \tilde A_{kj}) - a_{jk}(A_{kj} A_ {ji} - \tilde A_{kj} \tilde A_{ji}) \big)  \\
&\hspace{0.5cm}+ \big( a_{jk} (A_{ki} - \tilde A_{ki}) - a_{ik} (A_{ji} A_{ki} - \tilde A_{ji} \tilde A_{ji} )        \big)  \\
&\hspace{0.5cm}+ \big(   a_{jk}(A_{ik} -\tilde A_{ik}) - a_{jk}(A_{ji} A_{ik} - \tilde A_{ji} \tilde A_{ik})     \big) \\
& = : \mathcal J_{3k,1} + \mathcal J_{3k,2} + \mathcal J_{3k,3} + \mathcal J_{3k,4}. 
\end{aligned}
\end{align}
Since $\mathcal J_{3k,j},~j=1,\cdots,4$ has the same form, it suffices to  provide the estimate for  $\mathcal J_{3k,1}$:
\begin{align} \label{C-12-4}
\begin{aligned}
\mathcal J_{3k,1}& = a_{ik} ( A_{jk} - \tilde A_{jk}) - a_{jk}(A_{jk} A_{ji} - \tilde A_{jk} \tilde A_{ji}) \\
& = -a_{jk}(A_{ji} - \tilde A_{ji}) + a_{jk}(I_p - \tilde A_{jk} )(A_{ji} - \tilde A_{ji}) + a_{jk}(A_{jk} - \tilde A_{jk}) (I_p - A_{ji})  \\
&\hspace{0.5cm}+ (a_{ik} - a_{jk})(A_{jk} - \tilde A_{jk}).
\end{aligned}
\end{align}
 Thus, in \eqref{C-12-2}, we multiply $(A_{ji} - \tilde A_{ji})^\top $ and take the trace to find
\begin{align} \label{C-13}
\begin{aligned}
 &\frac12\frac{\d}{\d t} \|A_{ji} - \tilde A_{ji}\|_\tF^2  \\
 &\hspace{0.5cm}= \textup{tr} \big[ \{ (A_{ji} - \tilde A_{ji})\Xi_i - \Xi_j(A_{ji} - \tilde A_{ji} ) \}       (A_{ji} - \tilde A_{ji})^\top  \big] + \frac{\kp}{2N} \sum_{k=1}^N \textup{tr}[ \mathcal J_{3k}     (A_{ji} - \tilde A_{ji})^\top    ]  \\
 &\hspace{1cm}+ \frac{\kp}{2N} \sum_{k=1}^N a_{jk} \textup{tr}[ ( ( A_{ki} - \tilde A_{ki}) - (A_{ik} - \tilde A_{ik})    )   (A_{ji} - \tilde A_{ji})^\top ] \\
 &\hspace{1cm}+ \frac{\kp}{2N} \sum_{k=1}^N a_{ik}\textup{tr}[ ( ( A_{jk} - \tilde A_{jk}) - (A_{kj} - \tilde A_{kj})    )   (A_{ji} - \tilde A_{ji})^\top ]  \\
 &\hspace{0.5cm} = : \mathcal J_4 + \frac{\kp}{2N} \sum_{k=1}^N \mathcal J_{5k} + \frac{\kp}{2N} \sum_{k=1}^N \mathcal J_{6k} + \frac{\kp}{2N} \sum_{k=1}^N \mathcal J_{7k}.
\end{aligned}
\end{align}
Below, we estimate the terms ${\mathcal J}_4$ and $\mathcal J_{jk},~j=5,6,7$, separately. \newline

\noindent $\bullet$ (Estimate on $\mathcal J_4$): for a skew-symmetric matrix $Y$ and a matrix $B$, we observe the following identity:
\begin{equation*}
\textup{tr}[YBB^\top ] =0.
\end{equation*}
Then, we use the skew-symmetricity of $\Xi_i$  to see 
\begin{align*}
\mathcal J_4 &= \textup{tr} \big[      \{      (A_{ji} - \tilde A_{ji})\Xi_i - \Xi_j(A_{ji} - \tilde A_{ji} )  \}       (A_{ji} - \tilde A_{ji})^\top  \big] \\
& = \textup{tr} \big[   \Xi_i (    A_{ji} - \tilde A_{ji})^\top  (A_{ji} - \tilde A_{ji})   \big] -  \textup{tr} \big[   \Xi_j    (A_{ji} - \tilde A_{ji})(A_{ji} - \tilde A_{ji})^\top      \big] \\
& = 0.
\end{align*}

\noindent $\bullet$ (Estimate on $\mathcal J_{5k}$): Due to   \eqref{C-12-3} and \eqref{C-12-4}, it suffices to estimate the term involving $\mathcal J_{3k,1}$. Note that 
\begin{align*}
\textup{tr}[  \mathcal J_{3k,1}(A_{ji} - \tilde A_{ji})^\top  ]  & \leq  -a_{jk} \|A_{ji} - \tilde A_{ji}\|^2 + a_{jk} \sqrt p \big[ D(\mathcal S) + D( \tilde{\mathcal S}) \big]  d(\mathcal S, \tilde{\mathcal S})^2 \\
&\hspace{0.5cm}+ (a_{ik} - a_{jk}) \textup{tr} (( A_{jk} - \tilde A_{jk})(A_{ji} - \tilde A_{ji})^\top ),
\end{align*}
where we used the inequality:
\begin{equation*}
\|I_p - A_{ji} \|_\tF = \|S_j^\top S_j - S_j^\top  S_i\|_\tF \leq \sqrt p D(\mathcal S).
\end{equation*}
Hence, $\mathcal J_{5k}$ can be estimated as follows:
\begin{align*}
\mathcal J_{5k}&  \leq  \textup{tr}[ (\mathcal J_{3k,1}  +\cdots +  \mathcal J_{3k,4})(A_{ji} - \tilde A_{ji})^\top ]  \\
&=-2 (a_{ik} + a_{jk})\|A_{ji} - \tilde A_{ji}\|^2  + 2(a_{jk} + a_{ik} )\sqrt p [ D(S) + D(\tilde S) ] d(S,\tilde S)^2 \\
& \hspace{0.3cm}+ (a_{ik} - a_{jk}) \textup{tr} (( A_{jk} - \tilde A_{jk})(A_{ji} - \tilde A_{ji})^\top )  + (a_{ik} - a_{jk}) \textup{tr} (( A_{kj} - \tilde A_{kj})(A_{ji} - \tilde A_{ji})^\top ) \\ 
&\hspace{0.3cm} + (a_{jk} - a_{ik}) \textup{tr} (( A_{ki} - \tilde A_{ki})(A_{ji} - \tilde A_{ji})^\top )  + (a_{ij} - a_{ik}) \textup{tr} (( A_{ik} - \tilde A_{ik})(A_{ji} - \tilde A_{ji})^\top )  \\
& \leq -2 \kp a_m \|A_{ji} -\tilde A_{ji}\|_\tF^2 + 2\kp a_M \sqrt p \big[ D(\mathcal S) + D( \tilde{\mathcal S}) \big] d(\mathcal S,\tilde{\mathcal S})^2 + \mathcal J_{5k,1},
\end{align*}
where sum of the last four terms are denoted by $\mathcal J_{5k,1}$. Then, we observe
\begin{align*}
&\mathcal J_{5k,1} + \mathcal J_{6k} + \mathcal J_{7k}\\
 &\hspace{0.3cm} = a_{jk} \textup{tr}[ ( ( A_{ki} - \tilde A_{ki}) - (A_{ik} - \tilde A_{ik})    )   (A_{ji} - \tilde A_{ji})^\top ] \\
 &\hspace{0.8cm}+ a_{ik}\textup{tr}[ ( ( A_{jk} - \tilde A_{jk}) - (A_{kj} - \tilde A_{kj})    )   (A_{ji} - \tilde A_{ji})^\top ] \\
&\hspace{0.8cm}+ (a_{ik} - a_{jk}) \textup{tr} (( A_{jk} - \tilde A_{jk})(A_{ji} - \tilde A_{ji})^\top )  + (a_{ik} - a_{jk}) \textup{tr} (( A_{kj} - \tilde A_{kj})(A_{ji} - \tilde A_{ji})^\top ) \\ 
&\hspace{0.8cm} + (a_{jk} - a_{ik}) \textup{tr} (( A_{ki} - \tilde A_{ki})(A_{ji} - \tilde A_{ji})^\top )  + (a_{ij} - a_{ik}) \textup{tr} (( A_{ik} - \tilde A_{ik})(A_{ji} - \tilde A_{ji})^\top )   \\
&\hspace{0.3cm} \leq 4(N-1) ( a_M + d(\mathcal A))  d(\mathcal S,\tilde{\mathcal S})^2. 
\end{align*}
In \eqref{C-13}, we collect all estimates to obtain the desired inequality:
\begin{equation*}
 \frac{\d}{\d t} d(\mathcal S,\tilde{\mathcal S}) \leq -2\kp (\Lambda - a_M \sqrt p (D(S) + D(\tilde S)))   d(\mathcal S,\tilde{\mathcal S}), \quad \Lambda = \left( a_m - \frac{N-1}{N}(a_M + d(\mathcal A))\right) .
\end{equation*}

 \section{Proof of Lemma \ref{L4.2}} \label{sec:B}
 \setcounter{equation}{0}
In this appendix, we present a  proof of Lemma \ref{L4.2} in which the differential inequality for $\|\mathcal D_{ij}\|_\tF^2$ is derived. First, we observe the following second-order ODE of $\mathcal D_{ij}$:
\begin{align} \label{D-25-1}
\begin{aligned}
& m\ddot{\cD}_{ij} +\gamma\dot{\cD}_{ij} + \kp \xi_i \cD_{ij} + \kp ( \xi_i - \xi_j) S_j \\
& \hspace{0.5cm} = -m\left( S_i\dot{S}_i^\top \dot{S}_i -S_j\dot{S}_j^\top \dot{S}_j \right) +\left( S_i\Xi_i -S_j\Xi_j \right) \\
& \hspace{1cm} -\frac{m}{\gamma} \left( 2\dot S_i \Xi_i - S_i\Xi_i S_i^\top  \dot S_i + S_i \dot S_i^\top  S_i \Xi_i - 2\dot S_j \Xi_j + S_j \Xi_j S_j^\top  \dot S_j - S_j \dot S_j^\top   S_j \Xi_j  \right) \\
& \hspace{1cm} +\frac{\kappa}{N} \sum_{k=1}^N a_{jk} \left( \frac{1}{2}S_jS_j^\top S_k +\frac{1}{2}S_jS_k^\top S_j -S_j \right) -a_{ik} \left( \frac{1}{2}S_iS_i^\top S_k + \frac{1}{2}S_iS_k^\top S_i -S_i \right).
\end{aligned}
\end{align}
It follows from the definition of Stiefel manifold that 
\begin{align*}
\cD_{ik}^\top \cD_{ik} = (S_i- S_k)^\top  (S_i- S_k)= 2I_p -S_i^\top S_k -S_k^\top S_i,
\end{align*}
or equivalently,
\begin{equation} \label{D-25-2}
S_i -\frac{1}{2}S_iS_i^\top S_k -\frac{1}{2}S_iS_k^\top S_i = \frac{1}{2}S_i\cD_{ik}^\top \cD_{ik}.
\end{equation}
In \eqref{D-25-1}, we use the relation \eqref{D-25-2} to see
\begin{align}\label{D-27}
\begin{aligned}
& m\ddot{\cD}_{ij} +\gamma\dot{\cD}_{ij} +\kp \xi_i \cD_{ij} + \kp ( \xi_i - \xi_j) S_j \\
&\hspace{0.5cm} =-m( S_i\dot{S}_i^\top \dot{S}_i -S_j\dot{S}_j^\top \dot{S}_j ) +\left( S_i\Xi_i -S_j\Xi_j \right) +\frac{\kappa}{2N} \sum_{k=1}^N \left( a_{ik} S_i\cD_{ik}^\top \cD_{ik} -a_{jk} S_j\cD_{jk}^\top \cD_{jk} \right) \\
& \hspace{1cm} +\frac{m}{\gamma} \left( 2\dot S_i \Xi_i - S_i\Xi_i S_i^\top  \dot S_i + S_i \dot S_i^\top  S_i \Xi_i - 2\dot S_j \Xi_j + S_j \Xi_j S_j^\top  \dot Sj - S_j \dot S_j^\top   S_j \Xi_j  \right) \\
&\hspace{0.5cm}=: -m\cI_6+\cI_7 +\frac{\kappa}{2N} \sum_{k=1}^N \cI_{8k} +\frac{m}{\gamma}\cI_9.
\end{aligned}
\end{align}
By using \eqref{D-27}, we derive a second-order ODE of $\|\mathcal D_{ij}\|_\tF^2$:
\begin{align} \label{D-27-0} 
\begin{aligned}
& m\frac{\d^2}{\d t^2}\left\| \cD_{ij} \right\|_\tF^2 +\gamma\frac{\d}{\d t}\left\| \cD_{ij} \right\|_\tF^2 +2\kappa \xi_i \left\| \cD_{ij} \right\|_\tF^2 + 2\kp (\xi_i - \xi_j) \textup{tr}(I_p  -S_j^\top  S_i) \\
 &\hspace{0.3cm}= \textup{tr}\left[ m\left( \ddot{\cD}_{ij}^\top \cD_{ij} +2\dot{\cD}_{ij}^\top \dot{\cD}_{ij} +\cD_{ij}^\top \ddot{\cD}_{ij} \right) +\gamma\left( \dot{\cD}_{ij}^\top \cD_{ij} +\cD_{ij}^\top \dot{\cD}_{ij} \right) +2\kappa\cD_{ij}^\top \cD_{ij}\right] \\
& \hspace{0.3cm}= \textup{tr}\left[ \left( m\ddot{\cD}_{ij} +\gamma\ddot{\cD}_{ij} +\kappa\cD_{ij} \right)^\top  \cD_{ij} \right] +\textup{tr}\left[ \cD_{ij}^\top  \left( m\ddot{\cD}_{ij} +\gamma\ddot{\cD}_{ij} +\kappa\cD_{ij} \right) \right] +2m\| \dot{\cD}_{ij} \|_\tF^2 \\
  &\hspace{0.3cm}= 2\textup{tr} \left[ \cD_{ij}^\top  \left( -m\cI_6 +\cI_7 +\frac{\kappa}{2N} \sum_{k=1}^N \cI_{8k} +\frac{m}{\gamma}\cI_9 \right) \right] +2m\| \dot{\cD}_{ij} \|_\tF^2. 
\end{aligned} 
\end{align}
We now  assume
\begin{equation*}
\xi_i = \xi_j  \equiv \xi,\quad i,j=1,\cdots,N.
\end{equation*}
Then, we sum the relation \eqref{D-27-0} with respect to $i,j=1,\cdots,N$ to find temporal evolutions for $\displaystyle \mathcal G = \frac{1}{N^2} \sum_{i,j=1}^N \|\mathcal D_{ij}\|_\tF^2$:
\begin{align} \label{D-27-0-0}
\begin{aligned}
&m \ddot {\mathcal G} + \gamma \dot{\mathcal G} + 2\kp\xi  \mathcal G \\
&\hspace{0.5cm} =  \frac{2}{N^2} \sum_{i,j=1}^N \textup{tr}\left[   D_{ij}^\top  \left( -m\cI_6 +\cI_7 +\frac{\kappa}{2N} \sum_{k=1}^N \cI_{8k} +\frac{m}{\gamma}\cI_9 \right)  \right] + \frac{2m}{N^2} \sum_{i,j=1}^N \| \dot{\cD}_{ij} \|_\tF^2.
\end{aligned}
\end{align}
Below, we present estimates for $ \cI_k,~~ k=6,\cdots,9$, respectively. \newline

\noindent $\bullet$ (Estimate of $\cI_6$): we use the fact  
\begin{equation}  \label{D-27-1}
\|PAQ^\top \|_\tF = \|A\|_\tF\quad\textup{for}\quad P, Q \in \St\quad\textup{and}\quad A\in\cM_{p,p}(\bbr),
\end{equation}
which can be proved as follows:
\begin{align*}
\|PAQ^\top \|_\tF^2 = \textup{tr}( QA^\top  P^\top  PAQ^\top ) = \textup{tr}(QA^\top  AQ^\top  ) = \textup{tr}( A^\top AQ^\top Q) = \textup{tr}(A^\top A) = \|A\|_\tF^2.
\end{align*}
Then, we use \eqref{D-27-1} to  see
\begin{align*}
\begin{aligned}
\textup{tr} \left( \cD_{ij}^\top \cI_6 \right) &= \textup{tr} \big(( S_i^\top  -S_j^\top )( S_i\dot{S}_i^\top \dot{S}_i -S_j\dot{S}_j^\top \dot{S}_j ) \big) \\
& = \textup{tr} \big( \dot{S}_i^\top \dot{S}_i -S_i^\top S_j\dot{S}_j^\top \dot{S}_j -S_j^\top S_i\dot{S}_i^\top \dot{S}_i +\dot{S}_j^\top \dot{S}_j \big) \\
&= \| \dot{S}_i \|_\tF^2 +\| \dot{S}_j \|_\tF^2 -\textup{tr}\left( S_j\dot{S}_j^\top \dot{S}_jS_i^\top  +S_i\dot{S}_i^\top \dot{S}_iS_j^\top  \right) \\
&\leq 2( \| \dot{S}_i \|_\tF^2 +\| \dot{S}_j \|_\tF^2 ) \leq 4D ( \dot{\cS} )^2.
\end{aligned}
\end{align*}

\vspace{0.2cm}

\noindent $\bullet$ (Estimate of $\cI_7$): we use the maximality of $\|\Xi\|_\infty$ to find 
\begin{align*}
\begin{aligned}
\textup{tr} \left( \cD_{ij}^\top \cI_7 \right) &= \textup{tr} \left[ \left( S_i^\top -S_j^\top  \right)\left( S_i\Xi_i -S_j\Xi_j \right) \right] = \textup{tr} \left( \Xi_i -S_i^\top S_j\Xi_j -S_j^\top S_i\Xi_i +\Xi_j \right) \\
&\leq 2 \left\|\Xi_i\right\|_\tF +2\left\|\Xi_j\right\|_\tF \leq 4 \|\Xi\|_\infty.
\end{aligned}
\end{align*}

\vspace{0.2cm}

\noindent $\bullet$  (Estimate of $\cI_8$): by straightforward calculation, 
\begin{align*}
\frac{1}{N^3} \sum_{i,j,k}\textup{tr} ( \mathcal D_{ij}^\top  \mathcal I_{8k} ) &=\frac{1}{N^3} \sum_{i,j,k}a_{ik} \textup{tr} ( (S_i^\top  - S_j^\top )S_i \mathcal D_{ik}^\top  \mathcal D_{ik} )- a_{jk} \textup{tr}( (S_i^\top  - S_j^\top  ) S_j \mathcal D_{jk}^\top  \mathcal D_{ij}) \\
& = \frac{1}{2N^3} \sum_{i,jk} \Big( a_{ik} \textup{tr} (\mathcal D_{ik}^\top  \mathcal D_{ik}) - a_{jk} \textup{tr} (\mathcal D_{jk}^\top  \mathcal D_{jk})  \\
&\hspace{1cm}+ a_{jk} S_i^\top  S_j \mathcal D_{jk}^\top  \mathcal D_{jk} - a_{ik} S_j^\top  S_i \mathcal D_{ik}^\top  \mathcal D_{ik} \Big)= 0,
\end{align*}
where we used the interchange of the index $i \leftrightarrow j$.

\vspace{0.2cm}

\noindent $\bullet$  (Estimate of $ \cI_9$): we use $\|S\|_\tF = \sqrt{p}$ for $S\in\St$ to see
\begin{align*}
\textup{tr} \left( \cD_{ij}^\top \cI_9 \right) & = \textup{tr} \Big[ \left(S_i^\top  - S_j^\top \right)\big( 2\dot S_i \Xi_i - S_i \Xi_i S_i^\top  \dot S_i + S_i \dot S_i^\top  S_i \Xi_i  \\
&\hspace{2cm}-   2\dot S_j \Xi_j + S_j \Xi_j S_j^\top  \dot S_j - S_j \dot S_j^\top  S_j \Xi_j   \big)  \Big]  \\
& = \textup{tr}( S_i^\top  \dot S_i \Xi_i -\Xi_i S_i^\top  \dot S_i - 2S_i^\top  \dot S_j \Xi_j + S_i^\top  S_j \Xi_j S_j^\top  \dot S_j - S_i^\top  S_j \dot S_j^\top  S_j \Xi_j \\
& \hspace{1cm} + S_j^\top  \dot S_j \Xi_j - \Xi_j S_j^\top  \dot S_j -2S_j^\top  \dot S_i \Xi_i + S_j^\top  S_i \Xi_i S_i^\top  \dot S_i - S_j^\top  S_i \dot S_i^\top  S_i \Xi_i)  \\
& \leq 8 \sqrt p D(\dot{\mathcal S}) \|\Xi\|_\infty.
\end{align*}
%
In \eqref{D-27-0-0}, we combine all estimates  for  $\textup{tr} ( \cD_{ij}^\top \cI_k ),~k=6,\cdots,9$ and the fact that $\| \dot{\cD}_{ij} \|_\tF^2 \leq 4D(\dot{\cS})^2$ to obtain the desired differential inequality for $\mathcal G$:
\begin{align*}
m \ddot {\mathcal G} + \gamma \dot{\mathcal G} + 2\kp \xi \mathcal G  \leq 16mD(\dot{\cS})^2 +8\|\Xi\|_\infty    +\frac{16m\sqrt{p}\|\Xi\|_\infty}{\gamma}D(\dot{\cS}).
\end{align*}

\end{document}